\documentclass[11pt]{article}
\usepackage{amsthm}
\usepackage{amsfonts}
\usepackage[T1]{fontenc}
\usepackage{epsfig}

\usepackage[lflt]{floatflt}
\usepackage{epsfig}
\usepackage{graphicx}
\usepackage{amsmath}
\usepackage{amssymb}
\usepackage{color}
\usepackage{subfigure}
\usepackage{array}

\usepackage{geometry}
\usepackage{makeidx}
\usepackage{latexsym}

\newtheorem{note}{Note}
\newtheorem{theorem}{Theorem}
\newtheorem{definition}{Definition}
\newtheorem{proposition}{Proposition}

\newtheorem{remark}{Remark}

\newcommand{\dd}{{\rm d}}

\newcommand{\bbR}{{\mathbb R}}

\newcommand{\al}{\alpha}
\newcommand{\la}{\lambda}
\newcommand{\p}{\partial}
\newcommand{\be}{\beta}
\newcommand{\G}{\Gamma}

\newcommand{\DC}{_0^CD^\al_t}
\newcommand{\DRL}{_0^{RL}D^{1-\al}_t}

\begin{document}
\begin{center}\emph{}
\LARGE
\textbf{Explicit Solutions to Fractional Stefan-like problems for Caputo and Riemann--Liouville Derivatives}
\end{center}

                   \begin{center}
                  {\sc Sabrina D. Roscani$^{\sharp}${\footnote{This work was started at the beginning of 2018 when the second author was working also at Depto. Matem\'atica, FCEIA, UNR, Pellegrini 250,  Rosario, Argentina}} Nahuel D. Caruso$^{\dag,\ddag}$,  and Domingo A. Tarzia$^{\sharp}$}\\
									
$^{\sharp}$ CONICET - Depto. Matem\'atica,
FCE, Univ. Austral, Paraguay 1950, S2000FZF Rosario, Argentina \\
$^\dag$ Depto. Matem\'atica, EFB, UNR, Pellegrini 250,  Rosario, Argentina\\
$^\ddag$ CIFASIS - Centro Internacional Franco Argentino de Ciencias de la Informaci\'on y de Sistemas, CONICET, Bv. 27 de Febrero 210 Bis, Rosario, S2000EZP
Argentina\\

 (sroscani@austral.eud.ar, ncaruso@fceia.unr.edu.ar, dtarzia@austral.edu.ar)
                   \vspace{0.2cm}

       \end{center}
      
\small

\noindent \textbf{Abstract: }
  Two fractional  two-phase Stefan-like problems  are considered by using Riemann-Liouville and Caputo derivatives  of order $\al\in (0,1)$ verifying that they coincide with the same  classical Stefan problem at the limit case when $\al=1$. For both problems, explicit solutions in terms of the Wright functions are presented. Even though the similarity of the two solutions, a proof that  they are different is also given. The convergence when $\al \nearrow 1$  of the one and the other solutions to the same classical solution is given. Numerical examples for the dimensionless version of the problem are also presented and analyzed.



\section{Introduction}\label{Sec:Intro} 

This paper deals with Stefan--like problems governed by fractional diffusion equations (FDE). A classical Stefan problem is a problem  where a phase-change occurs, usually linked to  melting (change from solid to liquid) or freezing  (change from liquid to solid). In these problems the diffusion, considered as a heat flow, is expressed in terms of instantaneous local flow of temperature modeled by the Fourier Law. Therefore, the governing equations related to each phase are the well-known heat equations. There is also  a latent heat-type condition at the interface connecting the velocity of the free boundary and the heat flux of the temperatures in both phases known as ``Stefan condition''.
 A vast literature on Stefan problems is given in \cite{Alexiades, Cannon, Crank, Tarzia:biblio, Tarzia}. 
 
For example,  the following is the mathematical formulation for a classical one-dimensional two-phase Stefan problem: \textsl{Find the triple $\left\{ u_1,u_2, s\right\}$ such that they have sufficiently regularity and they verify that: 
}\begin{equation}{\label{St-Clasico}}
\begin{array}{llll}
     (i)  &   \frac{\p}{\p t}u_2(x,t)=\lambda^2_2\, \frac{\p^2}{\p x^2} u_2 (x,t), &   0<x<s(t), \,  0<t<T,  \, \,\\
     (ii)  &   \frac{\p}{\p t}u_1(x,t)=\lambda^2_1\, \frac{\p^2}{\p x^2} u_1 (x,t), &   x>s(t), \,  0<t<T,  \, \,\\
     (iii) &   u_1(x,0)=U_i, & 0\leq x,\\ 
       (iv)  &  u_2(0,t)=U_0,   &  0<t\leq T,  \\
         (v) & u_1(s(t),t)=u_2(s(t),t)=U_m, & 0<t\leq T, \\
  (vi) & \rho l \frac{d}{dt}s(t)=k_1\frac{\p}{\p x} u_1(s(t),t) -k_2\frac{\p}{\p x} u_2 (s(t),t), & 0<t\leq T, \\
	(vii) & s(0)=0,
                                             \end{array}
                                             \end{equation}
where $U_i<U_m<U_0$, $\lambda^2_j=\frac{k_j}{\rho c_j}$, $j=1$ (solid), $j=2$ (liquid) and we have assumed that  the thermophysical properties are constant as well as the free boundary can be represented by an increasing function of time.\\

Problem (\ref{St-Clasico}) is clearly governed by the heat equations $(\ref{St-Clasico}-i)$ and $(\ref{St-Clasico}-ii)$,  and 
has a phase-change condition (namely the Stefan condition) given by equation $(\ref{St-Clasico}-vi)$. \\

When the governing  equations $(\ref{St-Clasico}-i)$ and $(\ref{St-Clasico}-ii)$, or the Stefan condition $(\ref{St-Clasico}-vi)$ are replaced by other equations involving fractional derivatives in problems like $(\ref{St-Clasico})$,  we will refer to them as fractional Stefan-like problems. 

For example, the heat equation can be replaced by a fractional diffusion equation (FDE), which is closely linked to the study of anomalous diffusion.  A detailed explanation about the relation between anomalous diffusion and randon walk processes can be founded at the work done by Metzler and Klafter  \cite{MK:2000}. As we know, the diffusion equation is connected  to  the Brownian motion, where the mean square displacement (msd) of particles is  proportional to time. However, in Random Walks the msd is proportional to a power of  time. It is also interesting the approach given in \cite{BaFr:2005,GKS,Sax:1994} where it is suggested that anomalous diffusion could be caused by heterogeneities in the domain.

 For the relation between fractional diffusion equations and their applications,  we refer the reader to \cite{FM-Libro, Podlubny, Pskhu-Libro} and references therein where applications to the theory of linear viscoelasticity or thermoelasticity, among other, are presented. \\

In this paper, two approaches leading to subdiffusion are considered. The first one linked to 
the mathematical interest as generalized operators which interpolates classical derivatives 
(see \cite{Diethelm}), and  the second one related to Fourier's generalization laws (see \cite{Povstenko}).  These two approaches derived in two different formulations for the 
FDE. In order to present them, let a function $u=u(x,t)$ be  defined for given  one-dimensional variables $x$ and time $t$.
 A first formulation for the FDE  given in terms of  fractional integrals (see \cite{Fujita:1989}) is given by:
\begin{equation}\label{FDE-IntFrac} 
 _0I_t^\al u_{xx}(x,t) =  u(x,t) - u(x,0)
\end{equation}
where, $_0I^\al_t$ is the fractional integral of  Riemann--Liouville of order  {$\alpha$} in the $t-$variable defined as 
\begin{equation*}
_{0}I^{\alpha}_t u(x,t)=\frac{1}{\Gamma(\alpha)}\int^{t}_{0}(t-\tau)^{\alpha-1} u(x,\tau)\dd\tau 
\end{equation*}
for every $u$ such that $u(x,\cdot )\in L^1(0,T)$ for every $x>0$.
 Equation \eqref{FDE-IntFrac} is  derived also in \cite{MK:2000}, when a fractal time random walk is considered. As  it can be seen, no partial derivative in time is part of equation (\ref{FDE-IntFrac}),  but differenciating respect on time to both members we get a second formulation for a FDE
\begin{equation}\label{FDE-RL} 
  \,_0^{RL}D_t^{1-\al}u_{xx}(x,t)= u_{t}(x,t), 
\end{equation}
where $_0^{RL}D_t^{1-\al}$ is the fractional derivative of Riemann--Liouville in the $t-$variable defined  for every $\al \in (0,1)$ as
$$_{0}^{RL}D^{1-\al}_t u(x,t)=\frac{\p }{\p t}\, _{0}I^{\al}_t u(x,t) 
= \frac{1}{\Gamma(\al)}\frac{\p }{\p t}\int^{t}_{0}(t-\tau)^{\al-1} u(x,\tau)\dd\tau $$
for every $u\in AC_t[0,T]=\left\{ u \;|\;  u(x, \cdot ) \, \text{ is absolutely continuous on } 
[0,T] \text{ for every } x\in\bbR^+\right\}$.\\
Nevertheless,  when discussing about FDE associated to fractional time derivatives, 
the reader may retract on the FDE for the Caputo  derivative, that is
\begin{equation}\label{FDE-Caputo}
_{0}^CD^{\alpha}_t u(x,t)= u_{xx}(x,t). 
\end{equation}
Here, the partial time derivative has been replaced by a fractional derivative in the sense of Caputo respect on time. The Caputo derivative  $_0^{C}D_t^{\al}$ is defined for every $\al \in (0,1)$ as 
$$\,^C_{0} D^{\alpha}_t \,u(x,t)= \left[ \, _{0}I^{1-\alpha}_t \left( u_t \right)  \right] (x,t) = \frac{1}{\Gamma(1-\al)}\displaystyle\int^{t}_{0}(t-\tau)^{-\al} u_t(x,\tau)\dd\tau$$
 for every $u\in AC_t[0,T]$.

As we said before, in this paper, problems like (\ref{St-Clasico}) governed by equations like (\ref{FDE-RL}) or (\ref{FDE-Caputo}) will be studied. The literature on fractional phase-change problems is rather scant. In \cite{GeVoMiDa:2013} a fractional two-phase moving-boundary problem is approximated by a scale Brownnian motion model for subdiffusion. In \cite{Voller:2014} sharp and diffuse interface models of fractional Stefan problems are discussed. In \cite{RoBoTa:2018} a formulation of a one-phase fractional phase-change problem is given arising a time dependence on the initial extreme of the fractional derivative.  When the starting time considered in the fractional derivative of the governing equation is equal to 0, the mathematical point of view becomes interesting because they admit  self-similar solutions in terms of the Wright functions (see \cite{GoLuMa:2000, JiMi:2009, Pa:2013,RoSa:2013, RoTa:2014}). It is worth noting that this kind of problems are not deduced as in \cite{RoBoTa:2018, VoFaGa:2013}.\\  

This paper is a continuation of a previous work \cite{RoTa:2017-TwoDifferent}, related to  fractional   one-phase change problems. In Section \ref{Section:2} some basic definitions and properties on fractional calculus are given. In Section \ref{Section:3}, two fractional two-phase Stefan-like problems are considered, admitting both exact self-similar  solutions.  While the two governing equations are equivalent under certain assumptions for boundary-value-problems, when different ``fractional Stefan conditions'' are considered, the solutions obtained seem to be different. The uniqueness of the self-similar solution for one of the problems is obtained while it is an open problem for the other (see \cite{RoTa:2014}). Finally, numerical examples and graphics of the solutions are presented by 
considering a dimensionless model in Section \ref{Section:4}.

\section{Basic definitions and properties}
\label{Section:2}
\begin{proposition}\label{propo frac}\cite{Diethelm} The following properties involving the fractional integrals and derivatives hold:
\begin{enumerate}
\item \label{RL inv a izq de I} The  fractional derivative of  Riemann--Liouville  is a left inverse operator of the fractional integral of Riemann--Liouville of the same order  $\al\in \bbR^+$. If $f \in AC[a,b]$, then
$$^{RL}_{a}D^{\al}\,_{a}I^{\al}f(t)=f(t)  \quad \text{for every }\, t\in (a,b)$$

\item\label{not inverse} The fractional integral of Riemann--Liouville is not, in general, a left inverse operator of the fractional derivative of Riemann--Liouville.\\
 
In particular, if $0<\al<1$, then 
$ _{a}I^{\al}(^{RL}_{a}D^{\al}f)(t)=f(t) - \dfrac{_{a}I^{1-\al}f(a^+)}{\G(\al)(t-a)^{1-\al}}.$ 

\item\label{caso part I inv de RL} If there exist some $\phi \, \in L^1(a,b)$ such that $f=\,_aI^\al \phi$, then 
$$_{a}I^{\al}\,^{RL}_{a}D^{\al} f(t)=f(t)  \quad \text{for every } \, t\in (a,b).$$

\item\label{relacion RL-C} If  $f\in AC[a,b],$ then  
$$ ^{RL}_{a}D^{\al}f (t)=\frac{f(a)}{\G(1-\al)}(t-a)^{-\al}+\, ^C_{a} D^{\alpha}f(t).$$ 
\end{enumerate}
\end{proposition}

The fractional integral and derivatives of power functions can be easy calculated (see  e.g.  \cite{Podlubny}). In fact, for every $t\geq a$ we have that
\begin{equation}\label{Int-de-pot}
_aI^{\al}\left((t-a)^{\be}\right)=\frac{\G(\be + 1)}{\G(\be+\al+1)}(t-a)^{\be+\al }, \qquad \text{ for every } \quad  \be>-1, 
\end{equation}
and that
\begin{equation}\label{DRL-de-pot}
_a^{RL}D^{\al}\left((t-a)^{\be}\right)=\begin{cases}\frac{\G(\be + 1)}{\G(\be-\al+1)}(t-a)^{\be-\al} & \text{ if } \be \neq \al - 1,\\ 
0 & \text{ if } \be=\al-1.
\end{cases} 
\end{equation}
In particular, if $\be>0$, $_a^{RL}D^{\al}\left((t-a)^{\be}\right)=\,_a^{C}D^{\al}\left((t-a)^{\be}\right)$ due to Proposition \ref{propo frac} item \ref{relacion RL-C} and the Caputo derivative of $(t-a)^{\be}$ is not defined for $-1<\be<0$.

\begin{proposition} \cite{Samko} The following limits hold: 
\begin{enumerate}
\item If we set $_aI^0=Id$ for the identity operator, then for every $ f \, \in L^1(a,b)$,   
$$ \displaystyle\lim_{\al\searrow 0}\, _aI^\al f(t) =_aI^0f(t)=f(t), \qquad a.e.  $$
\item  For every $f \in AC[a,b]$, we have
$$ \displaystyle\lim_{\al\nearrow 1}\,_a^CD^\al f(t) = f'(t) \hspace{0.5cm} \text{and} \hspace{0.5cm}  \displaystyle\lim_{\al\searrow 1}\,_a^CD^\al f(t) = f'(t)-f'(a^+)\qquad \text{ for all } \, t\in (a,b). $$
\item  For every $f \in AC[a,b]$, 
$$ \displaystyle\lim_{\al\nearrow 1}\,_a^{RL}D^\al f(t) = f'(t)\hspace{0.5cm} \text{and} \hspace{0.5cm} \displaystyle\lim_{\al\searrow 1}\,_a^{RL}D^\al f(t) = f'(t) \quad a.e. $$
\end{enumerate}
\end{proposition}

\begin{definition} For every $x\in \bbR$ , the \textit{Wright} function is defined as
\begin{equation}\label{W} W(x;\rho;\be)=\sum^{\infty}_{k=0}\frac{x^{k}}{k!\G(\rho k+\be)} ,\quad  \rho>-1 \text{ and }  \be\in \bbR.\end{equation}
An important particular case of the Wright function is the \textit{Mainardi} function defined by 
$$ M_\rho (x)= W(-x,-\rho,1-\rho)=\sum^{\infty}_{n=0}\frac{(-x)^n}{n! \G\left( -\rho n+ 1-\rho \right)}, \quad \,0<\rho<1.  $$
\end{definition}

\begin{proposition}\label{Props W} \cite{Pskhu-Libro, Wr1:1934} 
  Let $\al > 0$, $\rho \in (0,1)$ and $\be \in \bbR$. Then the next assertions follows:
\begin{enumerate} 
\item \label{deriv comun} For every $x\in \bbR$ we have
$$ \frac{\p}{\p x} W(x,\rho,\be) = W(x,\rho,\rho+\be).$$  
\item For every $x>0$ and $c>0$,      
\begin{equation}\label{eq pskhu}
_0I^\al \left[ x^{\beta-1}W(-cx^{-\rho},-\rho, \beta)\right]=x^{\beta+\al-1}W(-cx^{-\rho},-\rho, \beta+\al) .
\end{equation}

\end{enumerate}
\end{proposition}
\begin{proposition}  \label{Propo monot Wrights}\cite{RoTa:2017-TwoDifferent,Wr1:1934} For every $\beta \geq 0$, $\rho \in (0,1)$: 
\begin{enumerate}
\item \label{W dec} The Wright function  $W\left(-\cdot,-\rho,\beta\right) $ is positive and strictly decreasing in $\bbR^+$.
\item \label{rel-Wrights} For every $x\geq 0$ the following equality holds  
$$\rho x W\left(-x,-\rho,\beta-\rho \right) = W\left(-x,-\rho,\beta- 1\right)+(1-\beta) W\left(-x,-\rho,\beta\right). $$ 
\item \label{desig-W} If, in addition  $0< \rho\leq \mu < \delta$, then for every $x>0$ the following inequality holds
\begin{equation}
\G(\delta)W(-x, -\rho, \delta)<\G(\mu) W(-x,-\rho, \mu).
\end{equation}

\end{enumerate}
\end{proposition}
%
            %

\begin{proposition}  \cite{Wr2:1940} \label{Wrights-x->infty}  For every $\beta \geq 0$ and  $\rho \in (0,1)$ the following limit holds
$$ \lim\limits_{x\rightarrow \infty}  W(-x,-\rho, \be)=0. $$
\end{proposition}
\begin{proposition}\label{convergencia-W}\cite{RoSa:2013, RoTa:2017-TwoDifferent}\label{conv M y W cuando al tiende a 1}Let $x\in \bbR^+_0$ be. Then  the following limits hold:\\
\begin{equation}\label{limite-M}\lim\limits_{\al\nearrow 1}M_{\al/2}\left(2x\right)=\lim\limits_{\al\nearrow 1}W\left(-2x,-\frac{\al}{2},1-\frac{\al}{2}\right)= M_{1/2}(2x)=\frac{e^{-x^2}}{\sqrt{\pi}},
\end{equation}

\begin{equation}\label{limite-W}
\lim\limits_{\al\nearrow 1}W\left(-2x,-\frac{\al}{2},\frac{\al}{2}\right)=\frac{e^{-x^2}}{\sqrt{\pi}},
\end{equation}
\begin{equation}\label{limite-erf}\lim\limits_{\al\nearrow 1}\left[1-W\left(-2x,-\frac{\al}{2}, 1\right)\right]={\rm erf}(x),
\end{equation} 
and
\begin{equation}\label{limite-erfc}\lim\limits_{\al\nearrow 1}\left[W\left(-2x,-\frac{\al}{2}, 1\right)\right]={\rm erfc}(x),
\end{equation}
where $erf(\cdot)$ is the error function defined by ${\rm erf}(x)=\frac{2}{\sqrt{\pi}}\displaystyle\int_0^xe^{-z^2}dz$ and ${\rm erfc}(\cdot)$ is the complementary error function defined by ${\rm erfc}(x)=1-{\rm erf}(x)$. Moreover, the convergence is uniform over compact sets.
\end{proposition}

\begin{proposition} The  fractional initial-boundary-value problems (\ref{Caputo-IC}) and (\ref{RL-IC}) for the quarter plane are equivalent if there exists $\be>0$ and $\delta >0$ such that  $\be<\al<1$  and $u_{xx}(x,\cdot)$ is an $O(t^{-\be})$ in $(0,\delta)$:   
\begin{equation}{\label{Caputo-IC}}
\begin{array}{llll}
     (i)  &  \DC u(x,t)= \frac{\p^2}{\p x^2} u (x,t), &   0<x, \,  0<t,  \, \,\\
     (ii)   & u(x,0)=u_0(x), & 0\leq x,\\ 
       (iii)  &  u(0,t)=g(t),   &  0<t,  
                                                      \end{array}
                                             \end{equation}

\begin{equation}{\label{RL-IC}}
\begin{array}{llll}
     (i)  &  \frac{\p}{\p t} u(x,t)=  \, \DRL \left( \frac{\p^2}{\p x^2} u (x,t)\right) , &   0<x, \,  0<t,  \, \,\\
     (ii)   & u(x,0)=u_0(x), & 0\leq x,\\ 
       (iii)  &  u(0,t)=g(t),   &  0<t,  
                                                      \end{array}
                                             \end{equation}
\end{proposition}
\proof Let $u=u(x,t)$ be a function satisfying equation $(\ref{Caputo-IC}-i)$. Applying $\DRL$ to both sides and using Proposition \ref{propo frac} item \ref{RL inv a izq de I} we get $(\ref{RL-IC}-i)$. \\
Let now, for the inverse suppose that $u$ satisfies equation $(\ref{RL-IC}-i)$. Applying $_0I^{1-\al}_t$ to both sides and using  Proposition \ref{propo frac} item \ref{not inverse} yields that 
 \begin{equation}\label{eq1}
 \DC u(x,t)= \frac{\p^2}{\p x^2} u (x,t)-
\frac{\displaystyle\lim\limits_{t\searrow 0}\, _0I^{\al}\left( \frac{\p^2}{\p x^2} u (x,t)\right)}{\Gamma(1-\al)t^{\al}}, \quad 0<x, \quad 0<t. 
\end{equation}
Now, for every $x$ fixed we have that  $u_{xx}(x,\cdot)$ is an $O(t^{-\be})$ in $(0,\delta)$, then for $t>0$ small it holds that 
\begin{equation}\label{Des-1} -C\tau^{-\be}\leq u_{xx}(x,\tau)\leq C\tau^{-\be},\quad 0<\tau\le t< \delta.  
\end{equation}
Multiplying by $\frac{(t-\tau)^{\al-1}}{\G(\al)}$ in (\ref{Des-1}),  integrating between 0 and $t$ and applying formula (\ref{Int-de-pot})  yields that   

\begin{equation}\label{Des-2} -C\frac{\G(1-\be) t^{\al-\be}}{\G(\al-\be + 1)}\leq \, _0I^{\al}_t \,u_{xx}(x,t) \leq C\frac{\G(1-\be) t^{\al-\be}}{\G(\al-\be + 1)},\quad \quad t<\delta.  
\end{equation}
Taking the limit when $ t$ tends to zero in $(\ref{Des-2})$ and being $\be<\al$ we conclude that 
equation $(\ref{Caputo-IC}-i)$ holds as we wanted to see.
\endproof

\begin{remark} Equations (\ref{Caputo-IC}-i) and (\ref{RL-IC}-i) has been treated as equivalent in 
literature, as it can be seeing at \cite{FM-Libro,   MK:2000, Povstenko}, but the condition 
\begin{equation}\label{cond FDE eq}
\displaystyle\lim\limits_{t\searrow 0}\, _0I^{\al}\left( \frac{\p^2}{\p x^2} u (x,t)\right)=0
\end{equation} 
must be considered and should not be forget it.
\end{remark}

\begin{remark}
It is easy to check that the following functions  verifies equation  (\ref{Caputo-IC}-i) and 
(\ref{RL-IC}-i) (we have taken $\lambda=1$ without loss of generality)
\begin{equation}\label{solFDE-1}
w_1(x,t)=x^2+ \frac{2}{\G(\al+1)}t^{\al}.
\end{equation}
\begin{equation}\label{solFDE-2}
w_2(x,t)=E_\al(t^{\al})\exp\left\{-x\right\}
\end{equation}
and 
\begin{equation}\label{solFDE-3}
w_3(x,t)=W\left(-\frac{x}{t^{\al/2}},-\frac{\al}{2}, 1\right).
\end{equation}
The condition (\ref{cond FDE eq}) trivially holds  for function $w_1$ and $w_2$ and it is no difficult to check it for  $w_3$ (by differenciating first and using Proposition \ref{Props W} then).
\end{remark}

\section{The Fractional Stefan-like Problems}
\label{Section:3}

In this section, two fractional Stefan-like problems admitting both explicit self-similar solutions will be treated. Before that, some clarification about the used terminology is presented. 

We refer to  fractional Stefan problems when the governed equations in such problem are derived  from physical assumptions, like considering memory fluxes.
 
For example, suppose that a process of  melting of a  semi--infinite slab  ($0\leq x<\infty$) of  some material  is taking place, and  the flux involved is a flux with memory. The melt temperature is $U_m$, and a constant temperature $U_0>U_m$ is imposed on the fixed face $x=0$. 
Let  $u_1=u_1(x,t)$ and $u_2=u_2(x,t)$ be the temperatures at the solid and liquid phases respectively. Let $J_1=J_1(x,t)$ and $J_2=J_2(x,t)$ be the respective functions for the fluxes at position $x$ and time $t$ and  let $x=s(t)$ be  the function representing the (unknown) position of the free boundary at time $t$. Suppose further that:\\
$(i)$ All the thermophysical parameters are constants.\\
$(ii)$ The function $s$ is an increasing function and consequently, an invertible function.\\
$(iii)$ $J_1 $ and $J_2$ are fluxes modeling the material with  memory which verifies that  ``\textsl{the  weighted sum of the fluxes back in time at the current time, is proportional to the gradient of temperature}'', that is, the following equations hold
\begin{equation}\label{J_1-elegido}
\nu_{\al}\,_{0}I^{1-\al}_tJ_1(x,t)=-k_1\frac{\p u_1}{\p x}(x,t)
\end{equation}
and 
\begin{equation}\label{J_2-elegido}
\nu_{\al}\,_{h(x)}I^{1-\al}_tJ_2(x,t)=-k_2\frac{\p u_2}{\p x}(x,t)
\end{equation}
where the initial time in the fractional integral (\ref{J_2-elegido}) is given by function  $h$ which gives us the time when the phase change occurs. That is,
$$ t=h(x)=s^{-1}(x) \qquad (\text{i.e. } \,\, x=s(t) )$$

 The number $\nu_{\al}$ is a parameter with physical dimension (see\eqref{med-nu_al}) such that 
\begin{equation}\label{nu-al}
\displaystyle\lim_{\al\nearrow 1}\nu_{\al}=1, \end{equation}  which has been added in order to preserve the consistency with respect to the units of measure in equations  (\ref{J_1-elegido}) and (\ref{J_2-elegido}). Also, the parameter 
\begin{equation}\label{mu-alpha} \mu_\al=\frac{1}{\nu_\al} 
\end{equation}
will be used in the following equations. More details about these parameters are given in Section 4.    \\

Making an analogous reasoning for the two-phase free-boundary problem, than the one made in \cite{RoBoTa:2018} for the one--phase free--boundary problem, the mathematical model for the problem described above is given by 

\begin{equation}{\label{Fr-Stefan-Problem}}
\begin{array}{llll}
     (i)  &   \frac{\p}{\p t}u_2(x,t)=\lambda_2^2\,\mu_{\al_2}\frac{\p}{\p x} \left( \,_{h(x)}^{RL}D^{1-\al}_t \left( \frac{\p}{\p x} u_2 (x,t)\right) \right), &   0<x<s(t), \,  0<t<T,  \, \,\\
     (ii)  &   \frac{\p}{\p t}u_1(x,t)=\lambda_1^2\,\mu_{\al_1}\frac{\p}{\p x} \left( _{0}^{RL}D^{1-\al}_t \left( \frac{\p}{\p x} u_1 (x,t)\right) \right), &   x>s(t), \,  0<t<T,  \, \,\\
     (iii) &   u_1(x,0)=U_i, & 0\leq x,\\ 
      (iv) &   u_2(0,t)=U_0,   &  0<t\leq T,  \\
         (v) & u_1(s(t),t)=u_2(s(t),t)=U_m, & 0<t\leq T, \\
  (vi) & \rho l\frac{d}{dt}s(t)= k_1\mu_{\al_1}\left. \DRL \frac{\p}{\p x} u_1(x,t)\right|_{(s(t)^+,t)}  & \\
  & \hspace{2.5cm} - k_2\mu_{\al_2}\left. _{h(x)}^{RL}D^{1-\al}_t \frac{\p}{\p x} u_2 (x,t)\right|_{(s(t)^-,t)}, & 0<t\leq T. \\
  s(0)=0
                                             \end{array}
                                             \end{equation}

where $U_i<U_m<U_0$ and $\mu_{\al}=\frac{1}{\nu_{\al}}$, (note that the parameter $\mu_\al$ can be the same in equations \eqref{Fr-Stefan-Problem}$-i$ and \eqref{Fr-Stefan-Problem}$-ii$,  
and without loss of generality we will take from now on that $\mu_{\al_2}=\mu_{\al_1}$).\\
Note that self-similar solutions to problem (\ref{Fr-Stefan-Problem}) had not been yet founded, due to the difficulty imposed by the variable button limit in the fractional derivative for the liquid phaace. 
As it was said at the beginning of this section, this paper deals with Stefan-like problems admitting explicit self-similar solutions. These problems come from the assumption of consider the button limit $t_0=0$ in the fractional time derivatives in the Caputo or Riemann--Liouville sense.    \\

\textbf{The Stefan-Like Problem for the Caputo derivative. }The next problem was treated in \cite{RoTa:2014} and can be obtained by replacing all the times derivatives in (\ref{St-Clasico}) by fractional derivatives in the Caputo sense of order $\al \in (0,1)$, i.e. 

\begin{equation}{\label{St-like-2ph-Caputo}}
\begin{array}{llll}
     (i)  &  \DC u_2(x,t)=\lambda^2_{\al_2}\, \frac{\p^2}{\p x^2} u_2 (x,t), &   0<x<s(t), \,  0<t<T,  \, \,\\
     (ii)  &   \DC u_1(x,t)=\lambda^2_{\al_1}\, \frac{\p^2}{\p x^2} u_1 (x,t), &   x>s(t), \,  0<t<T,  \, \,\\
     (iii) &   u_1(x,0)=U_i, & 0\leq x,\\ 
     (iv)  &   u_2(0,t)=U_0,   &  0<t\leq T,  \\
         (v) & u_1(s(t),t)=u_2(s(t),t)=U_m, & 0<t\leq T, \\
  (vi) & \rho l \DC s(t)=k_{\al_1} \frac{\p}{\p x} u_1(s(t)^+,t) - k_{\al_2}\frac{\p}{\p x} u_2 (s(t)^-,t), & 0<t\leq T, \\
(vii) & s(0)=0. & 	
                                             \end{array}
                                             \end{equation}

where $U_i<U_m<U_0$, $\lambda_{\al_i}$  are positive parameters named as ``subdiffusion coefficients'' given by 
$  \lambda_{\al_i}=\la_i\sqrt{\mu_{\al}}$ for  $i=1,2,$ 
and $ k_{\al_i}$ are positive parameters named as ``subdiffusion thermal conductivities'' given by $k_{\al_i}=k_i \mu_{\al}$, $i=1,2$.

\begin{definition}\label{Def sol St} The triple $\{u_1, u_2,s\} $ is a solution to  problem $(\ref{St-like-2ph-Caputo})$  if the following conditions are satisfied

\begin{enumerate}
    \item $u_1$  is continuous in the region $\mathcal{R_T}=\left\{(x,t)\colon 0\leq x \leq s(t), \, 0<t\leq T \right\}$ and at  the point $(0,0) $, $u_1 $ verifies that 
       $$ 0\leq \underset{(x,t)\rightarrow (0,0)}{\liminf}u_1(x,t)\leq \underset{(x,t)\rightarrow (0,0)}{\limsup } u_1(x,t)<+\infty .$$
			\item $u_2$  is continuous in the region $\left\{(x,t)\colon  x > s(t), \, 0<t\leq T \right\}$ and at  the point $(0,0) $, $u_2 $ verifies that 
       $$ 0\leq \underset{(x,t)\rightarrow (0,0)}{\liminf}u_2(x,t)\leq \underset{(x,t)\rightarrow (0,0)}{\limsup } u_2(x,t)<+\infty .$$

	\item  $u_1\in $ $C((0,\infty) \times(0,T) )\cap C^2_x((0,\infty) \times(0,T))$, such that $u_{1} \in \, AC_t[0,T]$
\item  $u_2\in $ $C((0,\infty) \times(0,T) )\cap C^2_x((0,\infty) \times(0,T))$, such that $u_{2} \in \, AC_t[0,T]$ 	.
  \item $s \in AC[0,T]$.
    \item $u_1$, $u_2$ and $s$ satisfy $(\ref{St-like-2ph-Caputo})$.
   	
		\end{enumerate}
\end{definition}

\begin{theorem}\cite{RoTa:2014} A self-similar solution to poblem (\ref{St-like-2ph-Caputo}) is given by 
\begin{equation}\label{sol Caputo-2f}
\left\{\begin{array}{l} u_2(x,t)=U_0-\frac{ U_0-U_m}{1-W\left(-2\xi_\al \lambda,-\frac{\al}{2},1\right)}\left[1-W\left(-\frac{x}{\lambda_{\al_2} t^{\al/2}},-\frac{\al}{2},1\right)\right] \\
        u_1(x,t)=U_i+\frac{U_m-U_i }{W\left(-2\xi_\al,-\frac{\al}{2},1\right)}W\left(-\frac{x}{\lambda_{\al_1} t^{\al/2}},-\frac{\al}{2},1\right) \\
				 s(t)=2\xi_\al  \lambda_{\al_1} t^{\al/2} \end{array}\right.
\end{equation}
 where $\xi_\al$ is \textbf{a} solution to the equation
\begin{equation}\label{eq xi}
\frac{k_{\al_2}(U_0-U_m)\G(1-\frac{\al}{2}) }{ \lambda_{\al_2}} F_2(2\lambda x)- \frac{k_{\al_1}(U_m-U_i) \G(1-\frac{\al}{2})}{\lambda_{\al_1}}F_1(2x)=\G\left(1+\frac{\al}{2}\right)\lambda_{\al_1} \rho l 2 x, \, x>0   \end{equation}

where  $ \lambda=\frac{\lambda_{\al_1}}{\lambda_{\al_2}}=\frac{\lambda_1\sqrt{\mu_\al}}{\lambda_2\sqrt{\mu_\al}} = \frac{\al_1}{\al_2}>0$,  and  $F_1:\bbR^+_0\rightarrow \bbR$ and $F_2:\bbR^+_0\rightarrow \bbR$  
are the functions defined by 
 \begin{equation}\label{F1-y-F2}
F_1(x)=\frac{M_{\al/2}(x)}{W\left(-x,-\frac{\al}{2},1\right)}  \quad \text{ and } \quad 
F_2(x)=\frac{M_{\al/2}(x)}{1-W\left(-x,-\frac{\al}{2},1\right)}  .
\end{equation}

\end{theorem}

%
%

\bigskip
\begin{note} The uniqueness of solution to equation $(\ref{eq xi})$ is still an open problem. However, the uniqueness of similarity solution will be achived next for the Riemann--Liouville Stefan--like problem. 
\end{note}

\noindent \textbf{The Stefan-Like Problem for the Riemann--Liouville derivative.} Consider now the following problem:
\begin{equation}{\label{St-like-2ph-RL}}
\begin{array}{llll}
     (i)  &   \frac{\p}{\p t}w_2(x,t)=\lambda^2_{\al_2}\,\frac{\p}{\p x} \left( \DRL \left( \frac{\p}{\p x} w_2 (x,t)\right) \right), &   0<x<r(t), \,  0<t<T,  \, \,\\
     (ii)  &   \frac{\p}{\p t}w_1(x,t)=\lambda^2_{\al_1}\,\frac{\p}{\p x} \left( \DRL \left( \frac{\p}{\p x} w_1 (x,t)\right) \right), &   x>r(t), \,  0<t<T,  \, \,\\
     (iii) &   w_1(x,0)=U_i, & 0\leq x,\\ 
       (iv)  &  w_2(0,t)=U_0,   &  0<t\leq T,  \\
         (v) & w_1(r(t),t)=w_2(r(t),t)=U_m, & 0<t\leq T, \\
  (vi) & \rho l \frac{d}{dt}r(t)= \left. k_{\al_1} \,\DRL \frac{\p}{\p x} w_1(x,t)\right|_{(r(t)^+,t)}  & \\
 &	\hspace{3cm} -k_{\al_2} \left. \DRL \frac{\p}{\p x} w_2 (x,t)\right|_{(r(t)^-,t)}, & 0<t\leq T,\\
			(vii)	& r(0)=0. &\\  
                                             \end{array}
                                             \end{equation}
where, as before,  $U_i<U_m<U_0$, $  \lambda_{\al_i}=\la_i\sqrt{\mu_{\al}}$ for  $i=1,2,$ and  $k_{\al_i}=k_i \mu_{\al}$, $i=1,2$.

\begin{remark}\label{cambio de limites en fr-St-Cond} The expression $\left. k_{\al_1} \,\DRL \frac{\p}{\p x} w_1(x,t)\right|_{(r(t)^+,t)} $ is equivalent to 
\begin{equation}\label{intercambio-lim-1} \lim\limits_{x\rightarrow r(t)^+} k_{\al_1}\,\DRL\, \frac{\p}{\p x} w_1(x,t),  \end{equation}
which should not coincide with 
\begin{equation}\label{intercambio-lim-2} k_{\al_1}\,\DRL\, \left( \lim\limits_{x\rightarrow r(t)^+} \frac{\p}{\p x} w_1(x,t)\right).  \end{equation}
\end{remark}

\begin{definition}\label{Def sol St} The triple $\{w_1, w_2,r\} $ is a solution of  problem $(\ref{St-like-2ph-RL})$  if the following conditions are satisfied

\begin{enumerate}
    \item $w_1$  is continuous in the region $\mathcal{R_T}=\left\{(x,t)\colon 0\leq x \leq s(t), \, 0<t\leq T \right\}$ and at  the point $(0,0) $, $u_1 $ verifies that 
       $$ 0\leq \underset{(x,t)\rightarrow (0,0)}{\liminf}w_1(x,t)\leq \underset{(x,t)\rightarrow (0,0)}{\limsup } w_1(x,t)<+\infty .$$
			\item $w_2$  is continuous in the region $\left\{(x,t)\colon  x > r(t), \, 0<t\leq T \right\}$ and at  the point $(0,0) $, $w_2 $ verifies that 
       $$ 0\leq \underset{(x,t)\rightarrow (0,0)}{\liminf}w_2(x,t)\leq \underset{(x,t)\rightarrow (0,0)}{\limsup } w_2(x,t)<+\infty .$$

	\item  $w_1\in $ $C((0,\infty) \times(0,T) )\cap C^2_x((0,\infty) \times(0,T))$, such that $w_{1x} \in 
	AC_t(0,T))$. 	
\item  $w_2\in $ $C((0,\infty) \times(0,T) )\cap C^2_x((0,\infty) \times(0,T))$, such that $w_{2x} \in 	AC_t[0,T]$ 	.
  \item $r \in C^1(0,T)$.
		\item There exist $\left.^{RL}_{0}D^{1-\al}_t \frac{\p}{\p x} w_2 (x,t)\right|_{(s(t)^+,t)}$ 
		and $\left.^{RL}_{0}D^{1-\al}_t \frac{\p}{\p x} w_1 (x,t)\right|_{(r(t)^-,t)}$ for all 
		$t \in (0,T]$.
    \item $w_1$, $w_2$ and $s$ satisfy $(\ref{St-like-2ph-RL})$.
   	
		\end{enumerate}
\end{definition}

\begin{theorem} An explicit solution for the two-phase fractional Stefan-like  problem (\ref{St-like-2ph-RL}) is given by
\begin{equation}\label{SOL-St-like-2ph-RL}
\left\{\begin{array}{l} w_2(x,t)=U_0-\frac{U_0-U_m}{1-W\left(-2\eta_\al \lambda,-\frac{\al}{2},1\right)}\left[1-W\left(-\frac{x}{\lambda_{\al_2} t^{\al/2}},-\frac{\al}{2},1\right)\right]\\
        w_1(x,t)=U_i+\frac{U_m-U_i}{W\left(-2\eta_\al,-\frac{\al}{2},1\right)}  W\left(-\frac{x}{\lambda_{\al_1} t^{\al/2}},-\frac{\al}{2},1\right) \\
				 r(t)=2\eta_\al  \lambda_{\al_1} t^{\al/2} \end{array}\right.
\end{equation}
 where $\eta_\al$ is the unique positive  solution to the equation
\begin{equation}\label{eq eta}
\frac{k_{\al_2}(U_0-U_m)}{ \lambda_{\al_1}\lambda_{\al_2}}G_2(2\lambda x) -\frac{k_{\al_1}(U_m-U_i)}{\lambda_{\al_1}^2 } G_1(2x)
 =\left(\rho l + \frac{k_{\al_1}(U_m-U_i)}{\lambda_{\al_1}^2} \right)2x,   
 \end{equation}
where  $ \lambda=\frac{\lambda_{\al_1\sqrt{\mu_\al}}}{\lambda_{\al_2}\sqrt{\mu_\al}}=\frac{\la_1}{\la_2}>0$, $U_i<U_m<U_0$  and  $G_1:\bbR^+_0\rightarrow \bbR$ and $G_2:\bbR^+_0\rightarrow \bbR$  
are the functions defined by 
 \begin{equation}\label{G1-y-G2}
G_1(x)=\frac{W\left(- x,-\frac{\al}{2},1+\frac{\al}{2}\right)}{W\left(- x,-\frac{\al}{2},1\right)}  \quad \text{ and } \quad      G_2(x)=\frac{2/\al W\left(-x,-\frac{\al}{2},\frac{\al}{2}\right)  }{1- W\left(-x,-\frac{\al}{2},1\right)}.
\end{equation}

\end{theorem}
\proof
Let the functions 

\begin{equation}{\label{sol-prop} }
\begin{array}{rcl}
w_i\colon  \bbR^+_0 \times (0,T)& \rightarrow & \bbR \\
(x,t) & \rightarrow  & w_i(x,t) = A_i + B_i \left[ 
1- W\left(-\frac{x}{\lambda_{\al_i} t^{\al/2}},-\frac{\al}{2}, 1 \right)\right]
\end{array}
\end{equation}

be the proposed solutions for $i=1,2.$ Rewriting expression  (\ref{eq pskhu}) for the variable $t$  and taking $c=\frac{x}{\lambda_{\al_i}}$  gives
\begin{equation}{\label{eq-pskhu-t}}
_0I^\al_t\,t^{\beta-1}W\left(-\frac{x}{\lambda_{\al_i}}t^{-\rho},-\rho, \beta\right)=t^{\beta+\al-1}W\left(-\frac{x}{\lambda_{\al_i}}t^{-\rho},-\rho, \beta+\al\right).
\end{equation}

Then, by using (\ref{eq-pskhu-t}) for $\be=1-\frac{\al}{2}$ and  Proposition \ref{Props W}  it is easy to check that $w_i$ verifies equations  $(\ref{St-like-2ph-RL}-i)$ and $(\ref{St-like-2ph-RL}-ii)$ respectively for $i=1,2.$\\

From condition $(\ref{St-like-2ph-RL}-v)$ we deduce that $r(t)$ must be proportional to $t^{\al/2}$. Therefore we set 
\begin{equation}\label{r}
r(t)=2\eta_\al \lambda_{\al_1} t^{\al/2}, \quad t \geq 0
\end{equation} 
where $\eta_\al$ is a constant to be determined and $\lambda_{\al_1}$ was added for simplicity in the next calculations. Now, from conditions $(\ref{St-like-2ph-RL}-iii)$, $(\ref{St-like-2ph-RL}-iv)$ and $(\ref{St-like-2ph-RL}-v)$ it holds that 
$$\begin{array}{ll}
 A_1= U_i+\frac{U_m-U_i}{W\left(-2\eta_\al,-\frac{\al}{2}, 1\right)}, 
& B_1= -\frac{U_m-U_i}{W\left(-2\eta_\al,-\frac{\al}{2}, 1\right)} \\
 A_2= U_0,  
& B_2= -\frac{U_0-U_m}{1-W\left(-2\eta_\al \lambda,-\frac{\al}{2}, 1\right)} 
\end{array} $$
As before, by considering (\ref{eq-pskhu-t}) for $\be=1-\frac{\al}{2}$ and Proposition \ref{Props W}, 
it holds that

$$\hspace{-10cm} \DRL w_{i_x}(x,t)= $$
\begin{equation}\label{DRL-w_x}
 \frac{B_i \al/2}{\lambda_{\al_1} \lambda_{\al_i}  t^{1-\al/2} } W\left(-\frac{x}{\lambda_{\al_i} t^{\al/2}},-\frac{\al}{2}, 1+\frac{\al}{2} \right)+ \frac{B_i \al/2}{\lambda_{\al_1} \lambda_{\al_i}}\frac{x}{t}W\left(-\frac{x}{\lambda_{\al_i} t^{\al/2}},-\frac{\al}{2}, 1 \right), \quad i=1,2.
\end{equation}

Then  replacing $(\ref{DRL-w_x})$ and $(\ref{r})$ in equation $(\ref{St-like-2ph-RL}-vii)$, and evaluating the limits following (\ref{intercambio-lim-1}) it yields that $\eta_\al $ must verify the next equality \\

\begin{equation}
\begin{split} \rho l 2 \eta_\al  \lambda_{\al_1}  = -\frac{k_{\al_1}(U_m-U_i)}{\lambda_{\al_1}^2}
     \frac{ W\left(-2\eta_\al,-\frac{\al}{2}, 1+\frac{\al}{2} \right)}{ W\left(-2\eta_\al,-\frac{\al}{2}, 1 \right)} - \frac{k_{\al_1}(U_m-U_i)}{\lambda_{\al_1}^2}2\eta_\al - \\
+\frac{k_{\al_2}(U_0-U_m)}{\lambda_{\al_1} \lambda_{\al_2}} \frac{ W\left(-2\lambda \eta_\al,-\frac{\al}{2}, 1+\frac{\al}{2} \right)}{1-  W\left(-\lambda 2\eta_\al,-\frac{\al}{2}, 1 \right)} + \frac{k_{\al_2}(U_0-U_m)}{\lambda_{\al_1}\lambda_{\al_2}}\frac{2 \lambda \eta_\al W\left(-\lambda 2\eta_\al,-\frac{\al}{2}, 1 \right)}{1-  W\left(-\lambda_\al 2\eta_\al,-\frac{\al}{2}, 1 \right)}. 
\end{split}
\end{equation} 
which leads to conclude that $\left\{w_1,w_2,r\right\}$ is a solution to (\ref{St-like-2ph-RL}) if and only if $\eta_\al $ is a solution to  the equation 

\begin{equation}\label{eq-eta-pre}
\begin{split}
\frac{k_{\al_2}(U_0 - U_m)}{\la_{\al_1} \la_{\al_2}}\frac{W\left(-\lambda 2 x,-\frac{\al}{2},1+\frac{\al}{2}\right) 
 + 2\la x W\left(-\lambda 2 x,-\frac{\al}{2},1\right)}{1-W\left(-\lambda 2 x,-\frac{\al}{2},1\right)}
 -\qquad \qquad\qquad\qquad
 \\
\qquad \qquad \qquad -k_{\al_1}\frac{U_m - U_i}{\la^2_{\al_1}}\frac{W\left(-2x,-\frac{\al}{2},1+\frac{\al}{2}\right)}{W\left(-2x,-\frac{\al}{2},1\right)}
=\left(\rho l +\frac{k_{\al_1}(U_m-U_i)}{\la_{\al_1}^2} \right)2x, \quad x>0.
\end{split}
\end{equation}

which, by using Proposition $\ref{Propo monot Wrights}-\ref{rel-Wrights}$ leads to equation (\ref{eq eta}).\\
The next step is to prove that Eq. (\ref{eq eta}) has unique solution. For that purpose 
we define function $G$ in $\bbR^+$ as
$$ G(x)= \frac{k_{\al_2}(U_0-U_m)}{ \lambda_{\al_1}\lambda_{\al_2}}G_2(2 \lambda x) -
          \frac{k_{\al_1}(U_m-U_i)}{\lambda_{\al_1}^2 } G_1(2x) -\left(\rho l+\frac{k_{\al_1}(U_m-U_i)}{\la_{\al_1}^2} \right)2x. $$

Note that $G$ is continuous function such that 
\begin{equation}\label{G-1a}
G(0^+)=+\infty. \end{equation}

From Proposition $\ref{Propo monot Wrights}-\ref{desig-W}$ for every $x>0$ we have that  
\begin{equation}\label{G2} 0<\frac{W\left(-2x,-\frac{\al}{2},1+\frac{\al}{2}\right)}{W\left(-2x,-\frac{\al}{2},1\right)}<\frac{1}{\G\left(\frac{\al}{2}+1\right)}, \end{equation}
then $G_1$ is bounded. Also, from (\ref{G2}) it holds that 
\begin{equation}\label{G-1b}
\begin{split} -\frac{k_{\al_1}(U_i-U_m)}{\lambda_{\al_1}^2 }\frac{1}{\G\left(\frac{\al}{2}+1\right)} + \frac{k_{\al_2}(U_0-U_m)}{ \lambda_{\al_1}\lambda_{\al_2}}G_2(2 \lambda x)-\left(\rho l+\frac{k_{\al_1}(U_m-U_i)}{\la_{\al_1}^2} \right)2x< \\
G(x)< \frac{k_{\al_2}(U_0-U_m)}{\lambda_{\al_1}\lambda_{\al_2}}G_2( 2\lambda x)-\left(\rho l+\frac{k_{\al_1}(U_m-U_i)}{\la_{\al_1}^2} \right)2x, 
\end{split}
\end{equation}
and taking the limit when $x \rightarrow \infty$ in (\ref{G-1b}) and using Proposition \ref{Wrights-x->infty} we obtain that 
\begin{equation}\label{G-1c}
G(+\infty)=-\infty. \end{equation}

Finally, consider the function $K\colon \bbR^+ \rightarrow \bbR$  defined as 
\begin{equation}\label{K}
K(x)=-\frac{k_{\al_1}(U_m-U_i)}{\lambda_{\al_1}^2 }\left[ G_1(2x) +2x\right] -\rho l 2x.
\end{equation} 

Applying Proposition \ref{Props W} item \ref{deriv comun} and being $\frac{(U_m-U_i)}{\lambda_{\al_1}^2 }>0$ it results that $K$ is a strictly decreasing function. By the other side, from Proposition \ref{Propo monot Wrights} item \ref{W dec}   we have that $G_2$ is a strictly decreasing function. Then it can be concluded that $G$ is a strictly decreasing function. Therefore Eq. (\ref{eq eta}) has a unique positive solution. 

\endproof
%

\begin{remark} The limits described in Remark \ref{cambio de limites en fr-St-Cond} are different if we compute them for the functions $w_1$ and $r$. In fact,  by using the computation made in the previous theorem,  we get\\

\begin{equation}\label{intercambio-lim-3}
\begin{split} \lim\limits_{x\rightarrow r(t)^+} \,\DRL\, \frac{\p}{\p x} w_1(x,t)=\frac{B_1}{\la_{\al_1}}\left[\frac{\al}{2}t^{\al/2-1}W\left(-2\eta_\al,-\frac{\al}{2},1+\frac{\al}{2}\right) + \right.\\
\left. \frac{\al}{2} 2\eta_\al t^{\al/2-1}W\left(-2\eta_\al,-\frac{\al}{2},1\right)  \right].  \end{split}\end{equation}

and from Proposition $\ref{Propo monot Wrights}-\ref{rel-Wrights}$, we have:  

\begin{equation}\label{Wr-al-2} W\left(-2\eta_\al,-\frac{\al}{2},1+\frac{\al}{2}\right)+ 2\eta_\al W\left(-2\eta_\al,-\frac{\al}{2},1\right) = \frac{2}{\al}  W\left(-2\eta_\al,-\frac{\al}{2},\frac{\al}{2}\right). \end{equation}

Then 
\begin{equation}\label{intercambio-lim-3'}
\lim\limits_{x\rightarrow r(t)^+} \,\DRL\, \frac{\p}{\p x} w_1(x,t)= \frac{B_1}{\la_{\al_1}} t^{\al/2-1}W\left(-2\eta_\al,-\frac{\al}{2},\frac{\al}{2}\right) 
\end{equation}

whereas

\begin{equation}\label{intercambio-lim-4}\, \DRL\, \left( \lim\limits_{x\rightarrow r(t)^+}\frac{\p}{\p x} w_1(x,t)\right)=\frac{B_1}{\la_{\al_1}} t^{\al/2-1} \frac{\G\left(1-\frac{\al}{2}\right)}{\G\left(\frac{\al}{2}\right)}M_{\al/2}(2\eta_\al).  \end{equation}

And we know that  $(\ref{intercambio-lim-3'})$ and $(\ref{intercambio-lim-4})$ are different due to  Proposition $\ref{Propo monot Wrights}-\ref{desig-W}$. 
\end{remark}

\begin{theorem}\label{differents} If $\la =1$,  the explicit solutions (\ref{SOL-St-like-2ph-RL}) to problem (\ref{St-like-2ph-RL}), and  (\ref{sol Caputo-2f}) to  problem (\ref{St-like-2ph-Caputo}) are different.
\end{theorem}
\proof Take $U_i=-1$, $U_m=0$ and $U_0=1$. Let $\left\{u_1, u_2, s\right\} $ be the solution to problem $(\ref{St-like-2ph-Caputo})$. Then $s(t)=2\la_{\al_1}\xi_\al t$ where   $\xi_\al$ is a  positive solution to equation  
\begin{equation}\label{eq xi-demos-dif-1}
\frac{k_{\al_2} \G(1-\frac{\al}{2}) }{\lambda_{\al_1} \lambda_{\al_2}} \frac{M_{\al/2}(2 \lambda x)}{1-W\left(-\lambda 2 x,-\frac{\al}{2},1\right)}- \frac{k_{\al_1} \G(1-\frac{\al}{2})}{\lambda^2_{\al_1}}\frac{M_{\al/2}(2x)}{W\left(-2x,-\frac{\al}{2},1\right)}=\G(1+\frac{\al}{2}) \rho l 2x.   
\end{equation}

By the other side, let $\left\{w_1,w_2,r\right\}$ be the solution to problem $(\ref{St-like-2ph-RL})$. Then $\eta_\al$ is the positive solution to equation

\begin{equation}
\frac{k_{\al_2} 2/\al}{ \lambda_{\al_1}\lambda_{\al_2}}\frac{ W\left(-2\lambda x,-\frac{\al}{2},\frac{\al}{2}\right) }{1- W\left(-\lambda 2x,-\frac{\al}{2},1\right)} -\frac{k_{\al_1}}{\lambda_{\al_1}^2 }\frac{W\left(- 2x,-\frac{\al}{2},1+\frac{\al}{2}\right)}{W\left(- 2x,-\frac{\al}{2},1\right)} =\left(\rho l + \frac{k_{\al_1}}{\lambda_{\al_1}^2} \right)2x,   
 \end{equation}

\noindent or equivalently, 

\begin{equation}\label{eq-RL-demos}
\begin{split}
\frac{k_{\al_2} }{ \lambda_{\al_1}\lambda_{\al_2}}\frac{ \G(1+\frac{\al}{2})2/\al  W\left(-2\lambda x,-\frac{\al}{2},\frac{\al}{2}\right) }{1- W\left(-2\lambda   x,-\frac{\al}{2},1\right)} -\qquad \qquad \qquad\qquad \qquad \qquad \qquad \qquad \qquad \qquad \qquad \qquad  \\ - \frac{k_{\al_1} \G(1+\frac{\al}{2})}{\lambda^2_{\al_1} }\frac{W\left(- 2x,-\frac{\al}{2},1+\frac{\al}{2}\right)+ 2x W\left(- 2x,-\frac{\al}{2},1\right)   }{W\left(- 2x,-\frac{\al}{2},1\right)}
 =\G(1+\frac{\al}{2})\rho l 2x.   
\end{split}
 \end{equation}

From Proposition $\ref{Propo monot Wrights}-\ref{rel-Wrights}$, for every $x>0$ we have that
\begin{equation}\label{caso part-igualdad W} W\left(- 2x,-\frac{\al}{2},1+\frac{\al}{2}\right)+ 2x W\left(- 2x,-\frac{\al}{2},1\right) = \frac{2}{\al}  W\left(- 2x,-\frac{\al}{2},\frac{\al}{2}\right). \end{equation}

Then using the fact that the Gamma function verifies that  $\frac{\G(1+\frac{\al}{2})}{\frac{\al}{2}}=\G(\frac{\al}{2})$ and replacing \eqref{caso part-igualdad W} in \eqref{eq-RL-demos} we deduce that  $\eta_\al$ is the unique positive solution to the equation 

\begin{equation}\label{eq nu-demos-dif-1-bis}
\frac{k_{\al_2}}{\la_{\al_1}\la_{\al_2}} \frac{ \G(\frac{\al}{2})  W\left(- 2\lambda  x,-\frac{\al}{2},\frac{\al}{2}\right) }{1- W\left(-2\lambda x,-\frac{\al}{2},1\right)} -\frac{k_{\al_1}}{\la^2_{\al_1}} \frac{ \G(\frac{\al}{2}) W\left(- 2x,-\frac{\al}{2},\frac{\al}{2}\right)}{W\left(-2 x,-\frac{\al}{2},1\right)}= \G(1+\frac{\al}{2})  \rho l  2x, \,x>0.   
 \end{equation}

If we suppose then that $\xi_\al=\eta_\al$, it result that there exist $\xi_\al >0$ such that 

\begin{equation}\label{eq-xi=eta}
\begin{split}
\frac{k_{\al_1}}{\la^2_{\al_1}}  \frac{ \G(\frac{\al}{2}) W\left(-2\xi_\al,-\frac{\al}{2},\frac{\al}{2}\right)}{W\left(-\xi_\al,-\frac{\al}{2},1\right)} - \frac{k_{\al_1}}{\la^2_{\al_1}} \frac{  \G(1-\frac{\al}{2}) M_{\al/2}(2\xi_\al)}{W\left(-\xi_\al,-\frac{\al}{2},1\right)}=\qquad \qquad\\
 \qquad \qquad=
\frac{k_{\al_2}}{\la_{\al_1}\la_{\al_2}}  \frac{ \G(\frac{\al}{2})  W\left(-\lambda 2\xi_\al,-\frac{\al}{2},\frac{\al}{2}\right) }{1- W\left(-\lambda \xi_\al,-\frac{\al}{2},1\right)} -  c_2  \frac{ \G(1-\frac{\al}{2}) M_{\al/2}( \lambda 2\xi_\al)}{1-W\left(-\lambda \xi_\al,-\frac{\al}{2},1\right)}.
\end{split}
\end{equation}

By using the hypothesis that $\la =1$, we conclude that 

\begin{equation} \label{eq h-demos-dif}
 \frac{\frac{k_{\al_1}}{\la^2_{\al_1}} }{W\left(-\xi_\al,-\frac{\al}{2},1\right)}=\frac{\frac{k_{\al_2}}{\la_{\al_1}\la_{\al_2}}}{1-W\left(-\lambda \xi_\al,-\frac{\al}{2},1\right)},
\end{equation} 
which leads to 
\begin{equation} \label{demos-dif-2}
W\left(-\xi_\al,-\frac{\al}{2},1\right)=\frac{1}{1+\frac{k_{\al_2}\la_{\al_2} }{k_{\al_1}\la_{\al_1}}}.
\end{equation}

Replacing $(\ref{demos-dif-2})$ in equation  $(\ref{eq xi-demos-dif-1})$ yields that  
$$  \rho l \lambda_{\al_1} 2\xi_\al=0 $$
which leads to $\xi_\al=0,$ contradicting the fact that $\xi_\al>0$.
\endproof

\begin{note} It is worth noting that an analogous proof for Theorem \ref{differents} but considering $\la\neq 1$ does not holds. In fact, if we define the function $h_\al\colon \bbR^+ \rightarrow \bbR$ as $$
h_\al(x)= \G\left(\frac{\al}{2}\right) W\left(-x,-\frac{\al}{2},\frac{\al}{2}\right)- \G\left(1-\frac{\al}{2}\right) M_{\al/2}( x)
$$ 

\begin{center}
\begin{figure}[h]
\includegraphics[width=0.5\columnwidth]{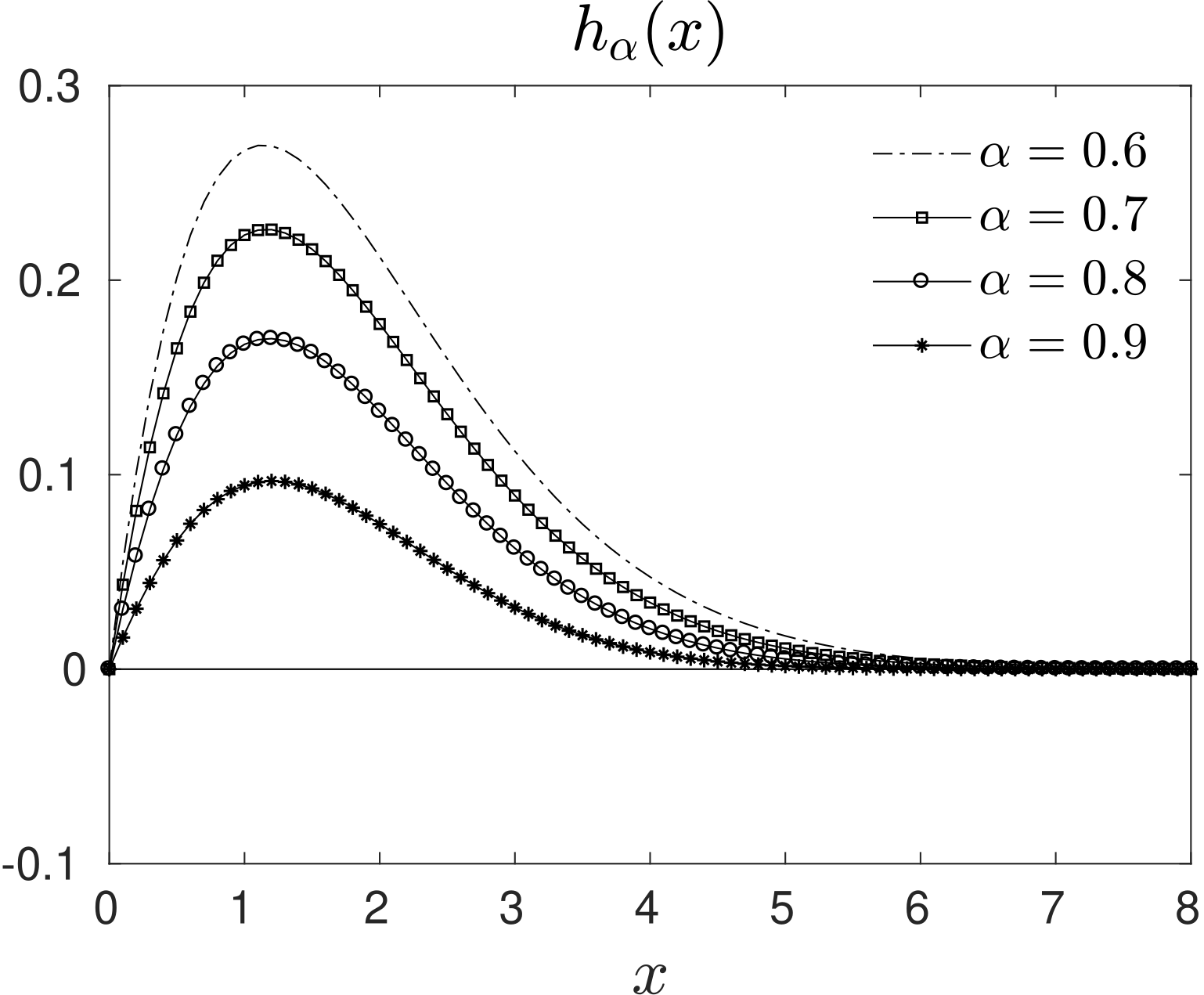} \\
\caption{The function $h_\al(x)= \G(\frac{\al}{2}) 
        W\left(-x,-\frac{\al}{2},\frac{\al}{2}\right)- \G(1-\frac{\al}{2}) M_{\al/2}( x)$ for different values of $\al$,}
\label{F:h-alpha}
\end{figure}
\end{center}

\noindent then equality (\ref{eq-xi=eta}) can be expressed as 
\begin{equation} \label{eq h-demos-dif}
\frac{k_{\al_2}}{\la_{\al_1}\la_{\al_2}}
 \frac{h_\al(\lambda 2\xi_\al)}{1-W\left(-\lambda 2\xi_\al,-\frac{\al}{2},1\right)}=
\frac{k_{\al_1}}{\la^2_{\al_1}}
  \frac{h_\al(2\xi_\al)}{W\left(- 2\xi_\al,-\frac{\al}{2},1\right)}.
\end{equation}
  If $\la \neq 1$, it is not possible to cancel the espressions $h_\al(\lambda 2\xi_\al)$ and 
  $h_\al(2\xi_\al)$ in equation (\ref{eq h-demos-dif}). Moreover the graphics in Figure \ref{F:h-quotient2-alpha} lead us to suppose that there exists a positive solution to equation     
	\begin{equation} \label{eq h-demos-dif-2}
\frac{k_{\al_2}}{\la_{\al_1}\la_{\al_2}} \frac{h_\al(\lambda x)}{1-W\left(-\lambda x,-\frac{\al}{2},1\right)}=
\frac{k_{\al_1}}{\la^2_{\al_1}}\frac{h_\al(x)}{W\left(- x,-\frac{\al}{2},1\right)}, \, x>0,
\end{equation}
then, it is not possible to get a contradiction like \eqref{demos-dif-2}.\\

%
%

\begin{center}
\begin{figure}[h]
\begin{tabular}{cc}
\includegraphics[width=0.5\columnwidth]{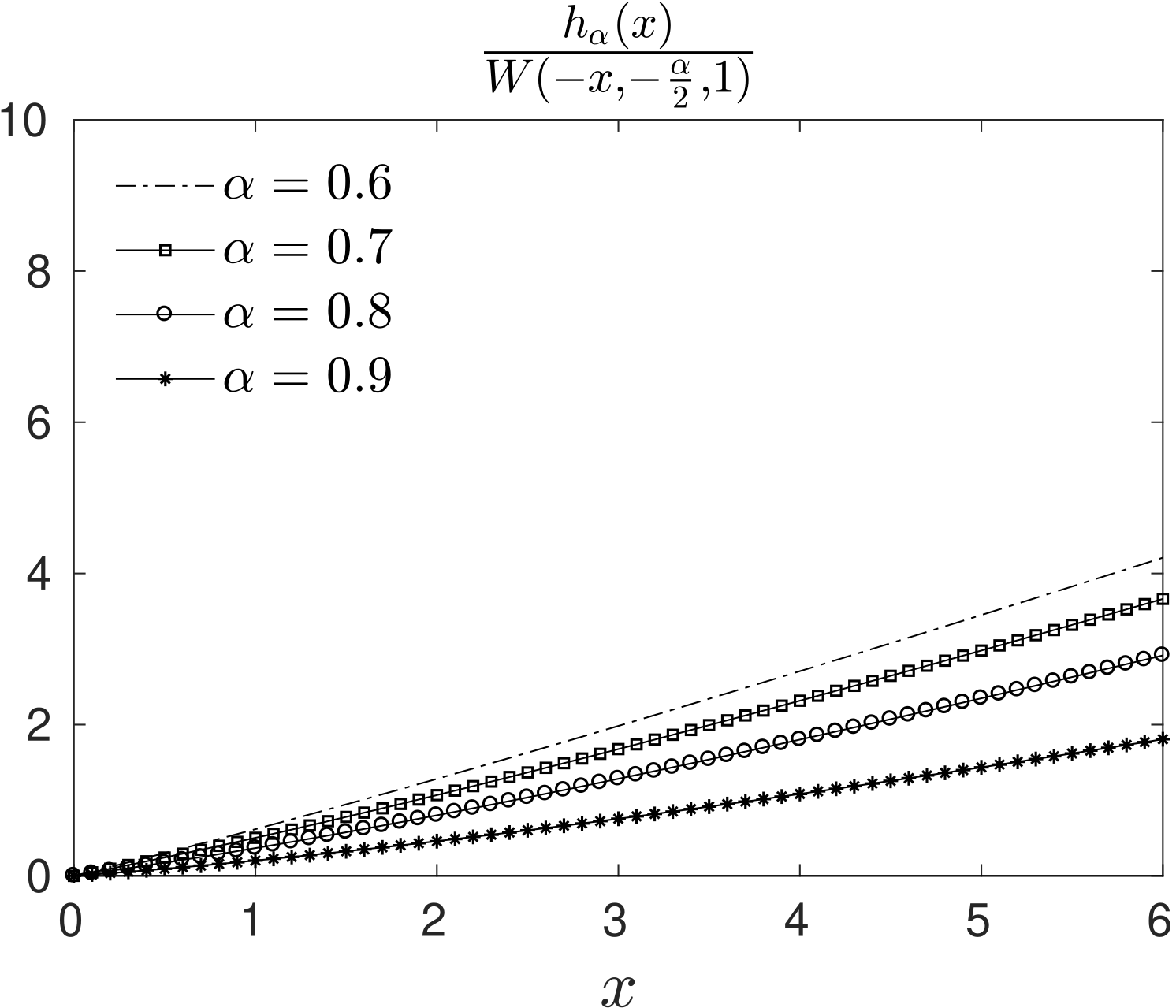} 
 &
\includegraphics[width=0.5\columnwidth]{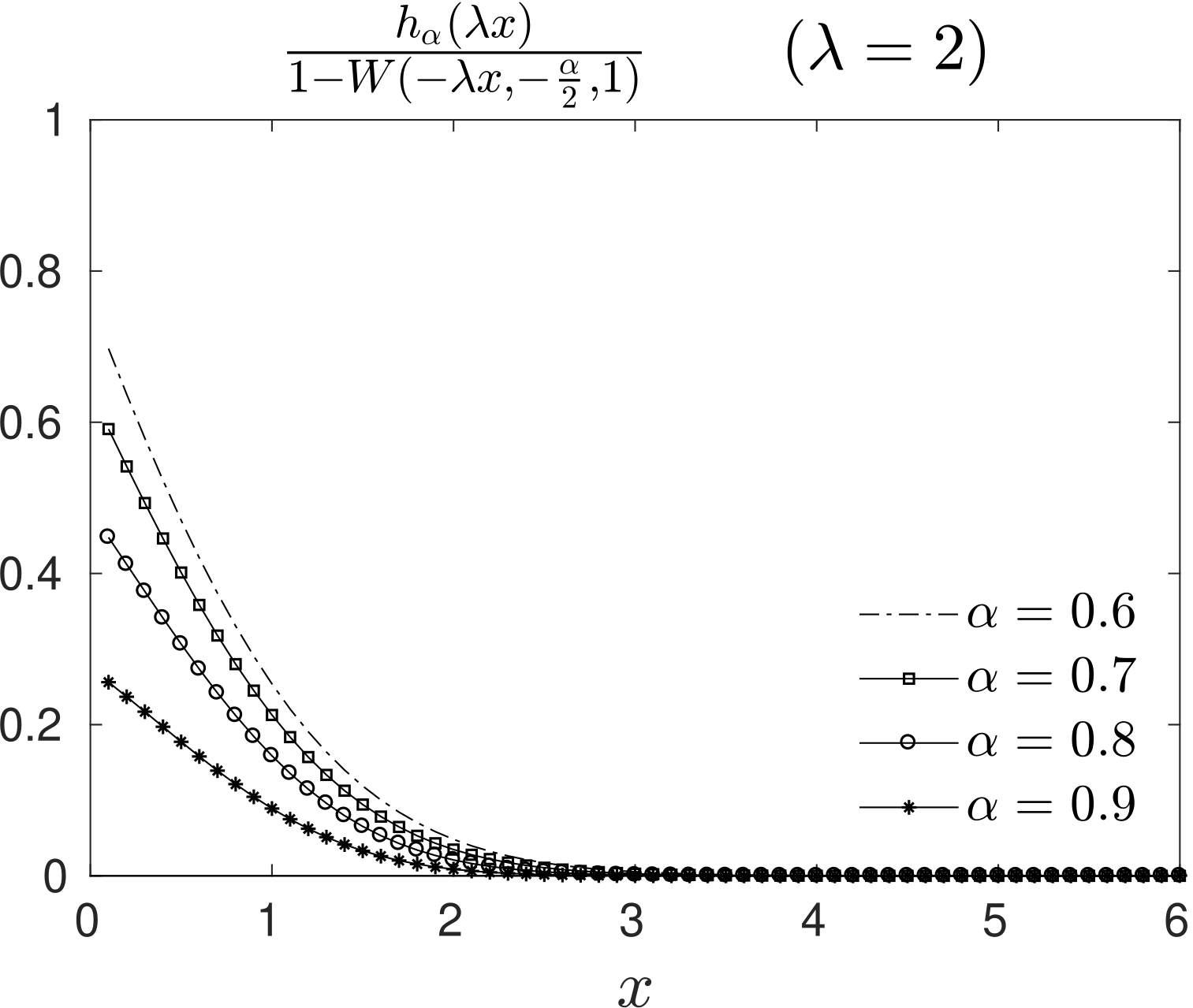} \\
\end{tabular}
\caption{The left and right quotients of equation (\ref{eq h-demos-dif-2})  for different values of $\al$}
\label{F:h-quotient2-alpha}
\end{figure}
\end{center}

However, if we take different values of $\la$ (which are different to 1) and the parameters $\xi_\al$ and $\eta_\al$ are estimated numerically for different values of $\alpha$, we show that they are different and converging both to the same value when $\al\nearrow 1$. Numerical examples will be given in the next section.
%
%

\end{note}

\begin{theorem} The explicit solution (\ref{SOL-St-like-2ph-RL}) to problem (\ref{St-like-2ph-RL})   converges, when $\al\nearrow 1$, to the unique solution to the classical Stefan problem given by
\begin{equation}{\label{Stefan-Clasico}}
\begin{array}{llll}
     (i)  &   \frac{\p}{\p t}u_2(x,t)=\lambda^2_2\, \frac{\p^2}{\p x^2} u_2 (x,t), &   0<x<s(t), \,  0<t<T,  \, \,\\
     (ii)  &   \frac{\p}{\p t}u_1(x,t)=\lambda^2_1\, \frac{\p^2}{\p x^2} u_1 (x,t), &   x>s(t), \,  0<t<T,  \, \,\\
     (iii) &   u_1(x,0)=U_i, & 0\leq x,\\ 
       (iv)  &  u_2(0,t)=U_0,   &  0<t\leq T,  \\
         (v) & u_1(s(t),t)=u_2(s(t),t)=U_m, & 0<t\leq T, \\
  (vi) & \frac{d}{dt}s(t)=k_1\frac{\p}{\p x} u_1(s(t),t) -k_2\frac{\p}{\p x} u_2 (s(t),t), & 0<t\leq T, \\
   (vii) & s(0)=0 
                                             \end{array}
                                             \end{equation}

\end{theorem}

\proof The unique solution to problem (\ref{Stefan-Clasico}) is the Neumann solution given in \cite{Weber:1901},

\begin{equation}
\left\{\begin{array}{l}
z_2(x,t)=U_0-(U_0-U_m)\frac{\mbox{erf\,}\left(\frac{x}{2\lambda_2\sqrt{t}}\right)}{\mbox{erf\,} \left(\nu_1\lambda\right)} \\
 z_1(x,t)=
     U_i+(U_m-U_i)\frac{\mbox{erfc\,} \left(\frac{x}{2\lambda_1\sqrt{t}}\right)}{\mbox{erfc\,} \left(\nu_1\right)} \\

w(t)= 2\eta_1 \lambda_1 \sqrt{t}
 \end{array}\right.
\end{equation}

\noindent where $ \eta_1$  is the unique  solution to the equation
  \begin{equation}\label{conv 4} \frac{k_2(U_0-U_m)}{ \lambda_1 \lambda_2}\frac{\exp\left\{ -\lambda^2 x^2 \right\}}{\sqrt{\pi} \mbox{erf\,}\left( \lambda x \right)}- \frac{k_1(U_m-U_i)}{\lambda_1^2}\frac{\exp \left\{-x^2\right\}}{\sqrt{\pi}\mbox{erfc\,} \left(x\right)} = \rho l x, \quad x>0. \end{equation}  

Reasoning like in  the previous theorem we can state that the solution  to problem (\ref{St-like-2ph-RL}) is given by (\ref{SOL-St-like-2ph-RL})  where $\eta_\al$ is the unique positive  solution to the equation\\

\begin{equation}\label{eq-eta-nueva}
\frac{k_{\al_2}(U_0-U_m)}{ \lambda_{\al_1}\lambda_{\al_2} \al}\frac{ W\left(-2\lambda  x,-\frac{\al}{2},\frac{\al}{2}  \right)}{1-W\left(-2\lambda  x,-\frac{\al}{2},1  \right)} -\frac{k_{\al_1}(U_m-U_i)}{\lambda_{\al_1}^2 \al } \frac{ W\left(-2 x,-\frac{\al}{2},\frac{\al}{2}  \right)}{W\left(-2x,-\frac{\al}{2},1  \right)}
 = \rho l x,   \quad x>0.
 \end{equation}

Clearly, if we take $\al=1$ in equation (\ref{eq-eta-nueva}) we recover equation (\ref{conv 4}). So, let the sequence  $\left\{ \eta_\al\right\}_\al$ be, where $\eta_\al$ is the unique positive solution to equation $(\ref{eq-eta-nueva})$ for each $0<\al<1$. \\
Defining the functions  
$$ f_\al(x)=\frac{k_{\al_2}(U_0-U_m)}{\rho l \lambda_{\al_1}\lambda_{\al_2} \al}\frac{ W\left(-2\lambda  x,-\frac{\al}{2},\frac{\al}{2}  \right)}{1-W\left(-2\lambda  x,-\frac{\al}{2},1  \right)} -\frac{k_{\al_1}(U_m-U_i)}{\rho l \lambda_{\al_1}^2 \al } \frac{ W\left(-2 x,-\frac{\al}{2},\frac{\al}{2}  \right)}{W\left(-2x,-\frac{\al}{2},1  \right)}
$$
for every $x \in \bbR^+$ and $0<\al \leq 1$, it holds that 
 $f_\alpha(\eta_\alpha)=\eta_\alpha$ for every $\alpha \in (0,1]$.  \\

From \cite{Tar:1981} we know that $f_1$ is a strictely decreassing function in $\bbR^+$. Taking a close interval $[a,b]\subset \bbR^+$ such that $\eta_1 \in [a,b]$, using the uniform convergence over compact sets  of all the positive functions given in Proposition \ref{convergencia-W}  and proceding like in \cite[Theorem 2]{RoTa:2017-TwoDifferent} we can state that
\begin{equation}\label{conv eta}
\lim\limits_{\al\nearrow 1}\eta_\al=\eta_1.
\end{equation}

Finally,   by taking the limit when $\al \nearrow 1$ in solution (\ref{SOL-St-like-2ph-RL}) by applying Proposition \ref{convergencia-W}, the thesis holds. 
\endproof

\begin{remark} By using the same technique described before, we can improve the result given in \cite[Theorem 3.3]{RoTa:2014} by considering the functions $g_\al$ defined  in $\bbR^+$ by
\begin{equation*}
\begin{split} g_\al(x)=\frac{k_{\al_2}(U_0-U_m)}{\rho l \lambda_{\al_1}\lambda_{\al_2} }\frac{\G(1-\al/2)}{\G(1+\al/2)}\frac{ M_{\al/2}\left(-2\lambda  x \right)}{1-W\left(-2\lambda  x,-\frac{\al}{2},1  \right)}\\
 -\frac{k_{\al_1}(U_m-U_i)}{\rho l \lambda_{\al_1}^2 \al }\frac{\G(1-\al/2)}{\G(1+\al/2)} \frac{ M_{\al/2}\left(-2 x\right)}{W\left(-2x,-\frac{\al}{2},1  \right)} 
\end{split}
\end{equation*}

and a sequence $\{\xi_\al\}_\al$  were $\xi_\al$ is a solution to 
equation $g_\al(x)=x$, $x>0$. 
\end{remark}

\section{The dimesionless problems and numerical results}
\label{Section:4}

In the aim to give different graphics of the solutions obtained in  Section 3, the  problems  (\ref{St-like-2ph-Caputo}) and (\ref{St-like-2ph-RL}) will be rewritten in their dimensionless form. \\

First, we give the following table exhibiting the usual heat conduction  physical dimensions related to this work. 
Let us write $\textbf{T}$ for temperature, $\textbf{t}$ for time, $\textbf{m}$ for mass and $\textbf{X}$ for position.

\begin{equation}\label{medidas}
\begin{array}{ccc}
u_1, u_2, w_1, w_2 & \textsl{temperatures} & [\textbf{T}]\\
k_1,k_2 & \textsl{thermal conductivities} & \left[\frac{{\textbf{m X}}}{{\textbf{Tt}^3}}\right] \\
\rho & \textsl{mass density} &  \left[\frac{{\textbf {m }}}{{\textbf {X}^3}}\right]  \\
c_1, c_2 & \textsl{specific heats} &   \left[\frac{\textbf{X}^2 }{\textbf{T}\textbf{t}^2}\right]\\
\la_i^2=\frac{k_i}{\rho c}, \, i=1,2 & \textsl{ diffusion coefficients} &  \left[\frac{\textbf{X}^2}{\textbf{t}}\right]\\
l & \textsl{latent heat per  unit  mass} & \left[\frac{\textbf{X}^2 }{\textbf{t}^2}\right]
\end{array}
\end{equation}

\begin{proposition}\label{frac-dim} For every $\al \in (0,1)$ it holds that
\begin{enumerate}
\item $\left[_0I^{\al} f \right]=[f]{\textbf{t}^\al}$ for every 
  $f=f(t) \in L^1(0,T)$.
\item $\left[_0^{RL}D^{\al} f \right]=\dfrac{[f]}{\textbf{t}^\al}$ for every $f=f(t) \in AC[0,T]$.
 \item   $\left[_0^{C}D^{\al} f \right]=\dfrac{[f]}{\textbf{t}^\al}$ for every $f=f(t) \in AC[0,T]$.
\end{enumerate}
\end{proposition}

Recall that the parameters $\nu_\al$ and $\mu_\al$ given in (\ref{nu-al})  where added  to preserve the consistency with respect to the units of measure in equations  (\ref{J_1-elegido}) and (\ref{J_2-elegido}). That is, being $ \left[ J \right]= \left[ k u_x \right]=\frac{\textbf{m}}{\textbf{t}^3}$ and using Proposition \ref{frac-dim}, it holds that 
 \begin{equation}\label{med-IalJ}
 \left[ _0I^{1-\al}_tJ(x,t) \right]=\left[\frac{1}{\G(1-\al)} \int_{0}^t\frac{J(x,\tau)}{(t-\tau)^\al}\dd\tau \right]=\frac{\textbf{m}}{\textbf{t}^{2+\al}}.  
 \end{equation}
 
Then, replacing (\ref{med-IalJ}) in (\ref{J_1-elegido}) one gets

\begin{equation}\label{med-nu_al}
\left[\nu_\al\right]=\frac{\left[k\frac{\p u}{\p x}\right]}{\left[_{h(x)}I^{1-\al}_tJ\right])}=\frac{1}{\textbf{t}^{1-\al}}.
\end{equation}

Therefore, \begin{equation}\label{med-mu_al}
\left[\mu_\al\right]=\textbf{t}^{1-\al}.
\end{equation}

\begin{proposition}\label{adim-fcns} Let $x_0$ be a characteristic position and let $U^*$ be a characteristic temperature. Then, if the following rescaling variable are considered 
\begin{equation}\label{variable-change}
y=\frac{x}{x_0}, \qquad \tau=\frac{\la_1^2}{x_0^2}t \quad  \text{and} \qquad \tilde{w}=\frac{w}{U^*},
\end{equation}
it holds that 
\begin{equation}\label{non-dim-1}
_0I^\al_t(w_{x}(x,t))=\frac{U^*x_0}{\la_1^2}\left(\frac{\la_1^2}{x_0^2}\right)^{1-\al}\,_0I^\al_\tau (\tilde{w}_{y}(y,\tau)),
\end{equation}
\begin{equation}\label{non-dim-2}
_0I^\al_t(w_{xx}(x,t))=\frac{U^*}{\la_1^2}\left(\frac{\la_1^2}{x_0^2}\right)^{1-\al}\,_0I^\al_\tau (\tilde{w}_{yy}(y,\tau))
\end{equation}
and 
\begin{equation}\label{non-dim-3}
_0^{RL}D^{1-\al}_t(w_{xx}(x,t))=\frac{U^*}{x_0^2}\left(\frac{\la_1^2}{x_0^2}\right)^{1-\al}\,_0^{RL}D^{1-\al}_\tau (\tilde{w}_{yy}(y,\tau)).
\end{equation}
\end{proposition}

\proof We prove here equation (\ref{non-dim-1}). By considering the rescaling (\ref{variable-change}), we have 
\begin{equation} \tilde{w}(y,\tau)=\frac{w(x(y),t(\tau))}{U^*}.\end{equation}

Then
$$ _0I^\al_t(w_{x}(x,t))= \frac{1}{\G(\al)}\int_0^t \frac{w_x(x,z)}{(t-z)^{1-\al}}dz= \frac{U^*}{\G(\al)}\int_0^t \frac{\frac{1}{x_0}\tilde{w}_y(y,\tau(z))}{(t-z)^{1-\al}}dz=$$
$$  = \frac{U^*}{\G(\al)}\int_0^{\frac{\la_1^2}{x_0^2}t} \frac{\tilde{w}_y(y,v)}{\left(\frac{x_0^2}{\la_1^2}\right)^{1-\al}(\frac{\la_1^2}{x_0^2}t-v)^{1-\al}}\frac{x_0}{\la_1^2}dv=  \frac{U^*}{x_0}\left(\frac{x_0^2}{\la_1^2}\right)^{\al}\,_0I^\al_\tau (\tilde{w}_{y}(y,\tau)).
  $$

\endproof

Now, let us consider problems problems  (\ref{St-like-2ph-Caputo}) and (\ref{St-like-2ph-RL}). By using Proposition \ref{adim-fcns} it is easy  to state that the  governing equation $(\ref{St-like-2ph-RL}-i)$ is equivalent to the following equation

\begin{equation}\label{equiv-2}
\frac{\p}{\p \tau}\tilde{w}_2(y,\tau)= \lambda^2 \mu_\al \left(\frac{\la_1^2}{x_0^2} \right)^{1-\al}\, _0^{RL}D^{1-\al}_\tau  \tilde{w}_{2yy} (y,\tau).
\end{equation} 

Note that $\mu_\al=\left(\frac{x_0^2 }{\la_1^2} \right)^{1-\al}$ is an admissible parameter because $[\mu_\al]=\textbf{t}^{1-\al}$ and that $\lim\limits_{\al\nearrow 1}\mu_\al=1$. 
Then, the parameter $\mu_\al \left(\frac{\la_1^2}{x_0^2} \right)^{1-\al}$ in equation (\ref{equiv-2}) can be omitted.  \\

Analogously, transforming the governing equations, the Stefan conditions and the initial and boundary data in problems (\ref{St-like-2ph-Caputo})  and (\ref{St-like-2ph-RL}), and by taking $U_m=0$ and  $U^*=|U_i|$,   it follows that  the non-dimensional   associated form are given by 

 \begin{equation}{\label{adim-St-like-2ph-Caputo}}
\begin{array}{llll}
     (i)  &  \,_0^{C}D^{\al}_\tau \tilde{u}_2(y,\tau)=\lambda^2 \tilde{u}_{2yy}(y,\tau), &   0<y<\tilde{s}(\tau), \,  0<\tau<\tilde{T},  \, \,\\
     (ii)  &   \,_0^{C}D^{\al}_\tau \tilde{u}_1(y,\tau)= \tilde{u}_{2yy}(y,\tau), &   y>\tilde{s}(\tau), \,  0<\tau<\tilde{T},  \, \,\\
     (iii) &   \tilde{u}_1(y,0)=-1, & 0\leq x,\\ 
     (iv)  &   \tilde{u}_2(0,\tau)=\frac{U_0}{|U_i|},   &  0<\tau\leq \tilde{T},  \\
         (v) & \tilde{u}_1(\tilde{s}(\tau),\tau)=\tilde{u}_1(\tilde{s}(\tau),\tau)=0, & 0<\tau\leq \tilde{T}, \\
  (vi) &  \,_0^{C}D^{\al}_\tau \tilde{s}(\tau)= {\rm{Ste}} \left[ \tilde{u}_{1y}(\tilde{s}(\tau)^+,\tau) - \frac{k_2}{k_1} \tilde{u}_{2y} (\tilde{s}(\tau)^-,\tau)\right], & 0<\tau\leq \tilde{T}, \\
(vii) & \tilde{s}(0)=0. & 	
                                             \end{array}
                                             \end{equation}
and
\begin{equation}{\label{adim-St-like-2ph-RL}}
\begin{array}{llll}
     (i)  &   \tilde{w}_{2\tau}(y,\tau)=\lambda^2 \,_0^{RL}D^{1-\al}_\tau {w}_{2yy} (y,\tau), &   0<y<\tilde{r}(\tau), \,  0<\tilde{t}<\tilde{T},  \, \,\\
     (ii)  &   \tilde{w}_{1\tau}(y,\tau)= \,_0^{RL}D^{1-\al}_\tau {w}_{1yy} (y,\tau), &   y>\tilde{r}(\tau), \,  0<\tau<\tilde{T},  \, \,\\
     (iii) &   \tilde{w}_1(y,0)=-1, & 0\leq y,\\ 
       (iv)  &  \tilde{w}_2(0,t)=\frac{U_0}{|U_i|},   &  0<\tau\leq \tilde{T},  \\
         (v) & \tilde{w}_1(\tilde{r}(\tau),\tau)=\tilde{w}_2(\tilde{r}(\tau),\tau)=0, & 0<\tau\leq \tilde{T}, \\
  (vi) &  \frac{d}{dt}\tilde{r}(\tau)= {\rm{Ste}} \left[ \left. \,_0^{RL}D^{1-\al}_\tau w_{1y}(y,\tau)\right|_{(\tilde{r}(\tau)^+,\tau)} \right.  & \\
 &	\hspace{3cm} \left. -\frac{k_2}{k_1} \left. \,_0^{RL}D^{1-\al}_\tau \tilde{w}_{2y} (y,\tau)\right|_{(\tilde{r}(\tau)^-,\tau)}\right], & 0<\tau\leq \tilde{T},\\
 			(vii)	& \tilde{r}(0)=0. &\\  

                                             \end{array}
                                             \end{equation}

where  $\la=\frac{\la_2}{\la_1}$ and ${\rm{Ste}}=\frac{|U_i|c_1}{l}$ is the non-dimensional Stefan number. 

\vspace{1cm}

In the following table there are different tests, i.e. sets of parameters for $\lambda$, 
$\frac{k_2}{k_1}$, $U=\frac{U_0}{|U_i|}$ and  $Ste$. For each test in Table 1 a correpondig graphic of the comparison  between the $\xi_{\al}$ and $\eta_{\al}$ is given in Figure 3.

\begin{table}[h!]
\centering
\begin{tabular}{|c |c | c| c| c|}
\hline
       & $\lambda$ & $\frac{k_2}{k_1}$ & $U=\frac{U_0}{|U_i|}$  &  $Ste$ \\
\hline
Test 1 & 0.5 & 0.5 & 1.0  &  0.5 \\
Test 2 & 2.0 & 2.0 & 1.0  &  0.5 \\
Test 3 & 0.5 & 0.5 & 1.0  &  1.2 \\
Test 4 & 2.0 & 2.0 & 1.0  &  1.2 \\
\hline
\end{tabular}
\caption{Different set of tests}
\end{table}

\newpage

\begin{figure}[h]
\label{xi_al-vs-eta_al}
\begin{tabular}{cc}
\includegraphics[width=0.45\columnwidth]{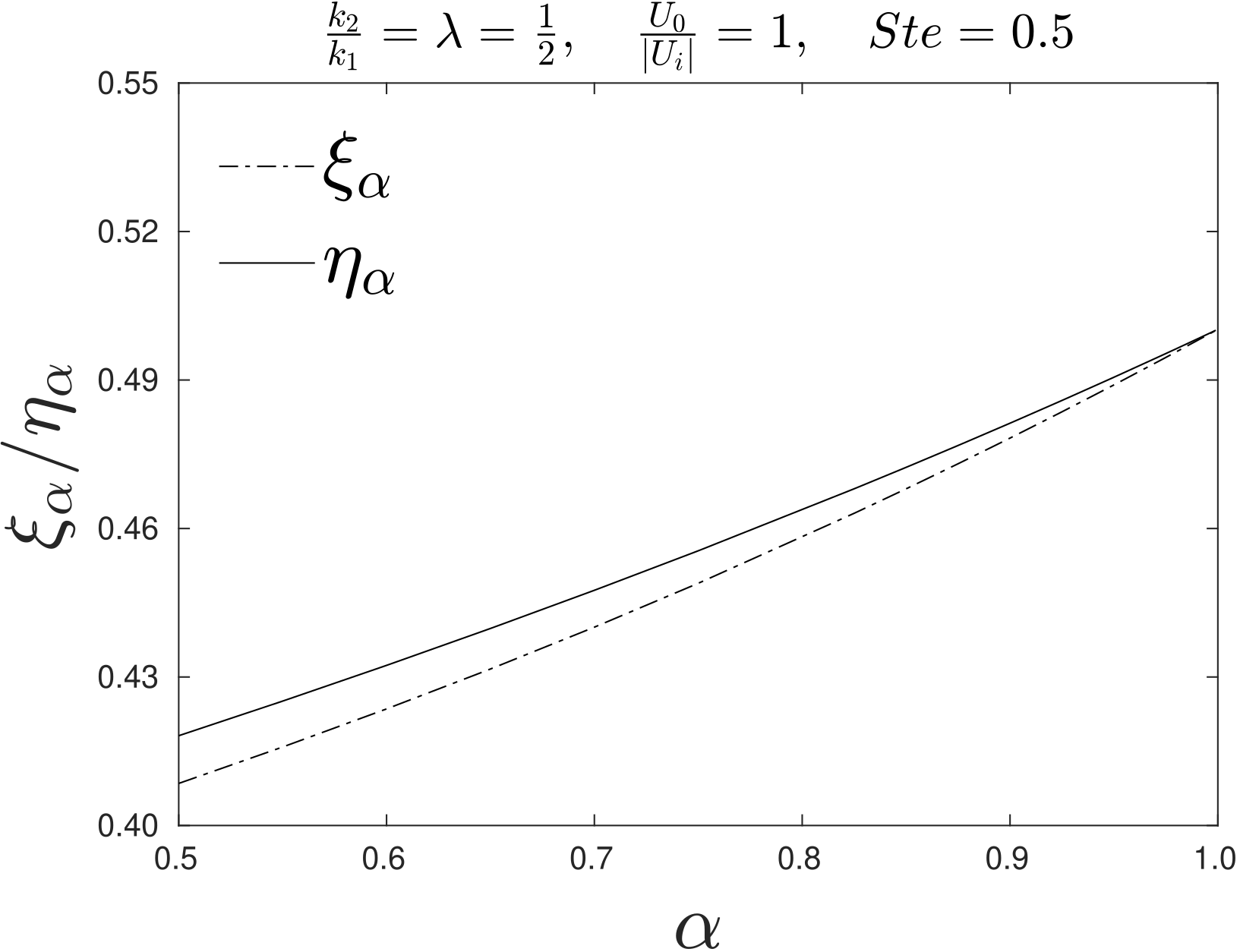} &
\includegraphics[width=0.45\columnwidth]{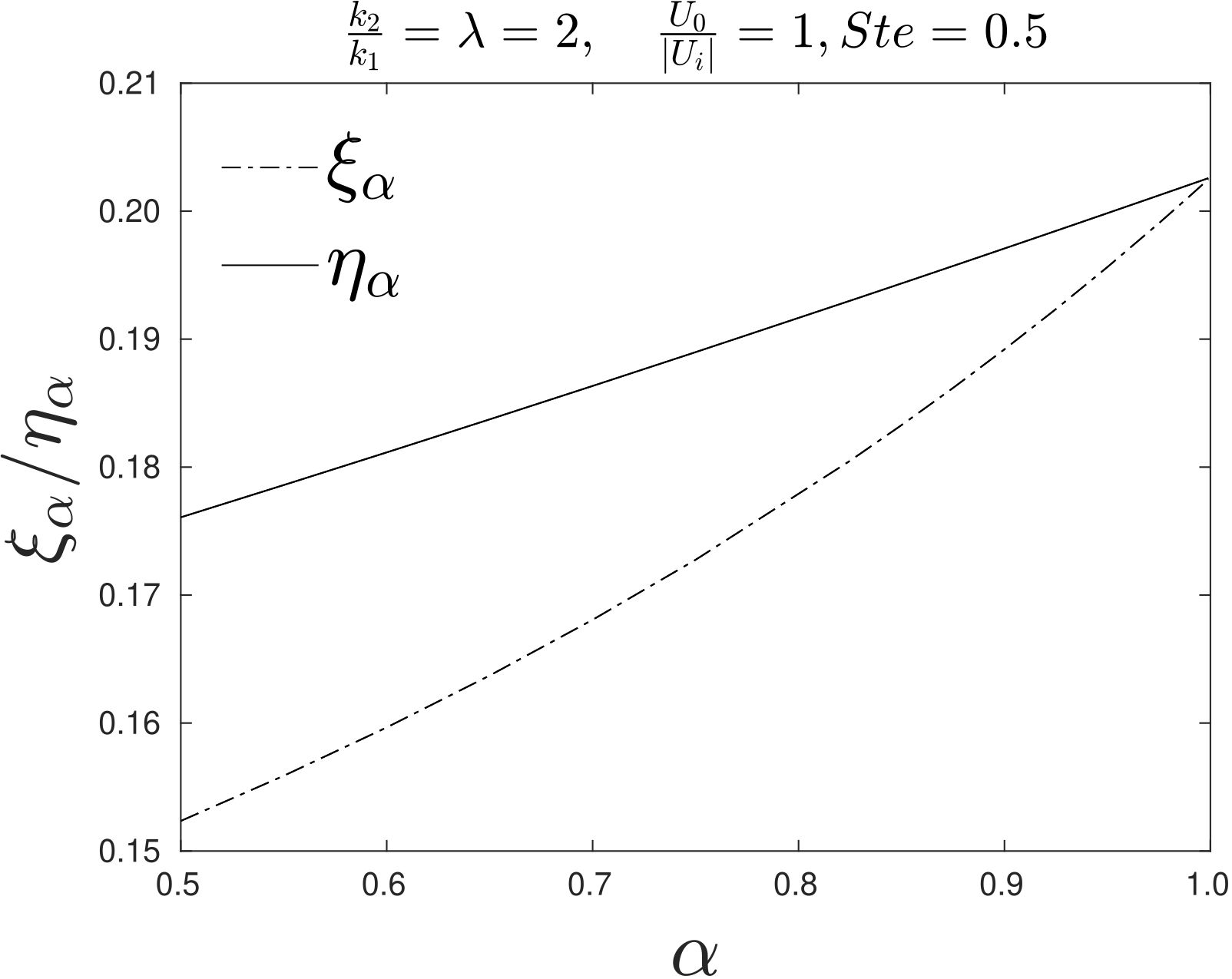} \\
\includegraphics[width=0.45\columnwidth]{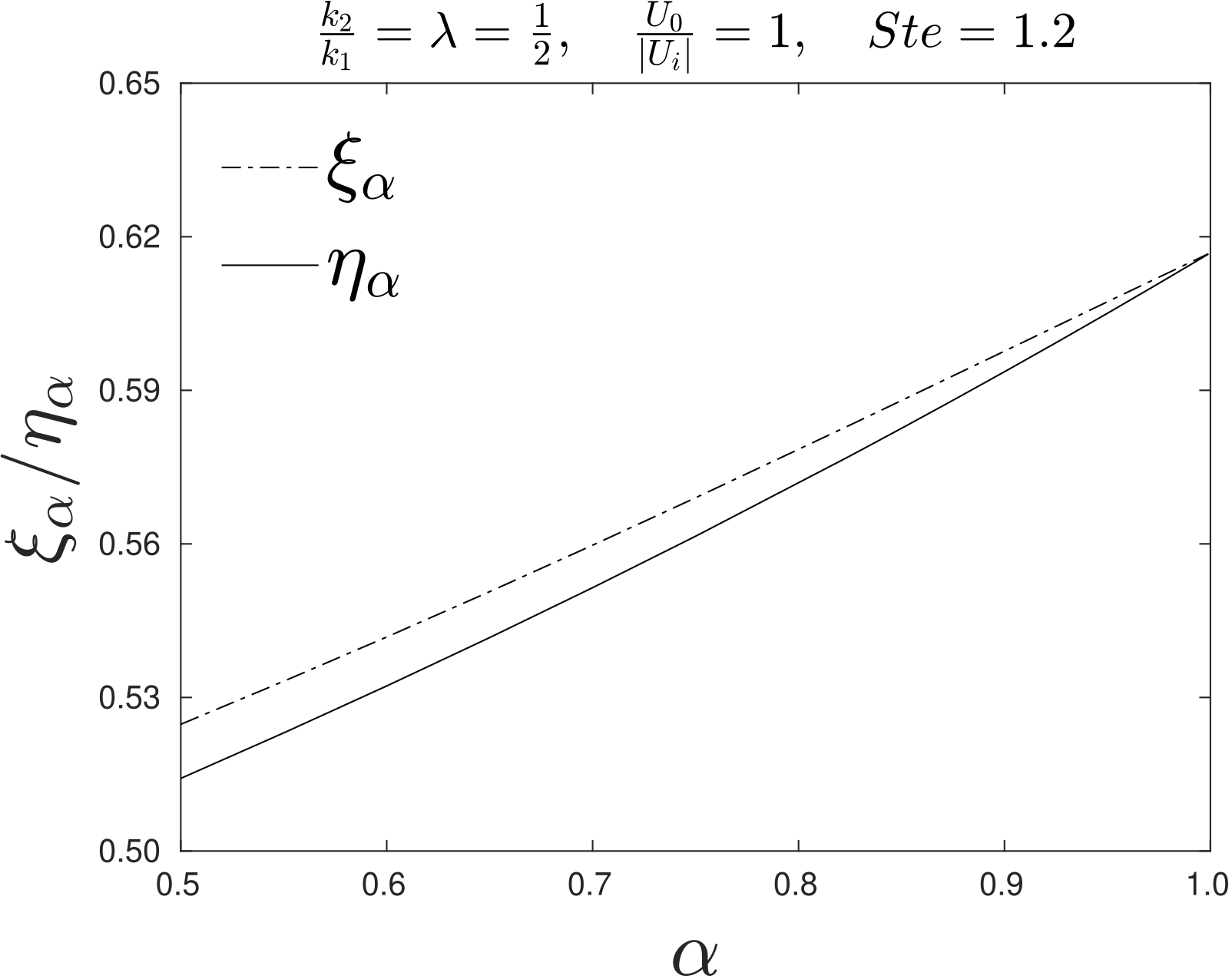} &
\includegraphics[width=0.45\columnwidth]{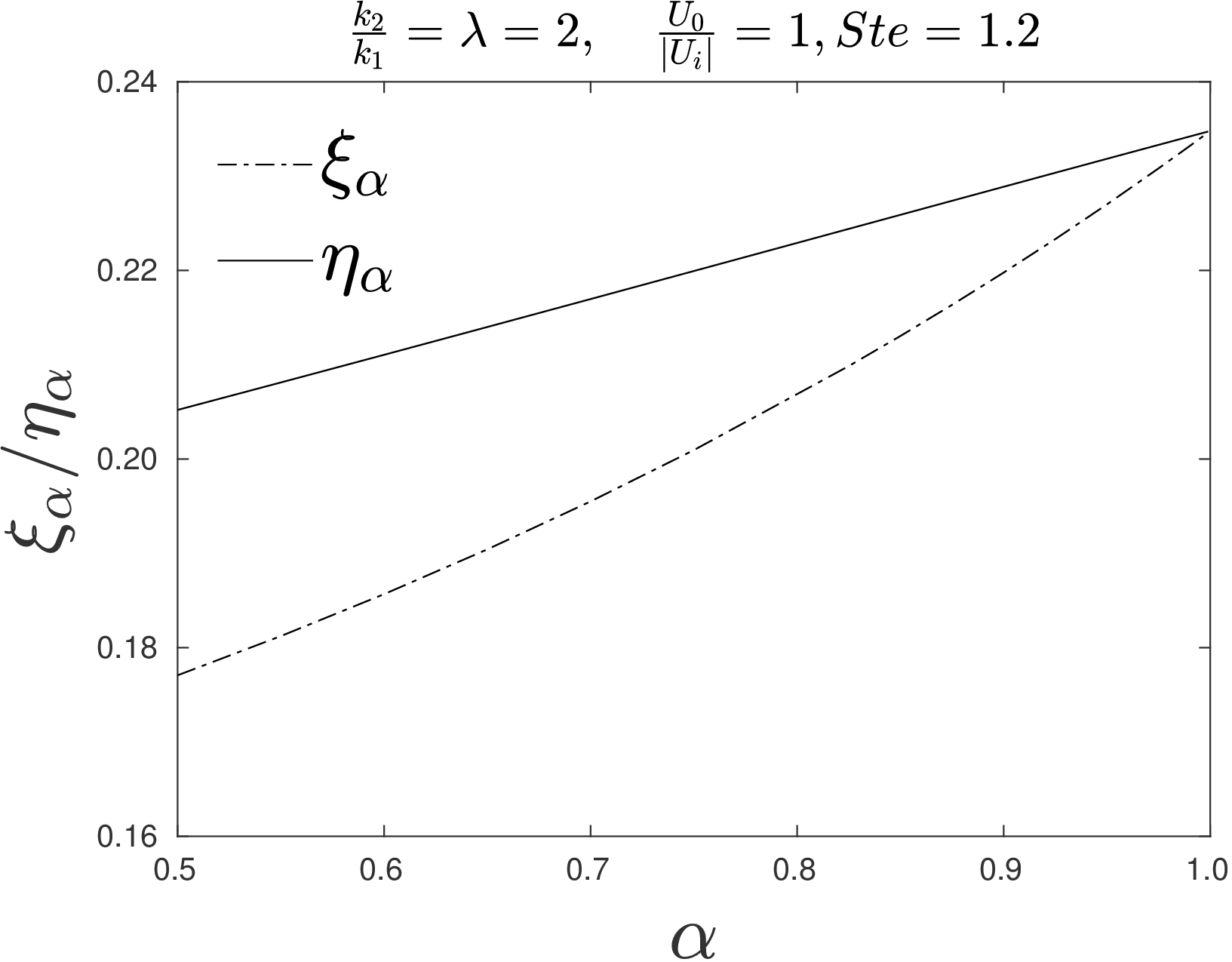} \\
\label{xi-vs-eta}
\end{tabular}
\caption{$\xi_\al$ vs. $\eta_\al$ for different values of $\al$}
\end{figure}

%

 At the end, we present  in Figures 4 and  5  some color maps of temperature   for  tests 2 and  3, respectively. Three values of $\al$ are considered and as it is expected from Theorem 4, both solutions approach themselves when $\al \nearrow 1$. 

\begin {figure}[h]
\label{C_vs_RL_test2}
\caption {Caputos's approach Solutions Vs. Riemann-Liouville's aproach Solutions for Test 2} 
\centering
\begin{tabular}{|m{0.5cm}|m{7cm}| m{7cm}|} 

\hline
 $\alpha$ & \textbf{Caputo Stefan-Like Pb. } & \textbf{Riemann-Liouville  Stefan-Like Pb. } \\
\hline   
\hline
$0.7$ & \includegraphics[width=0.5\columnwidth]{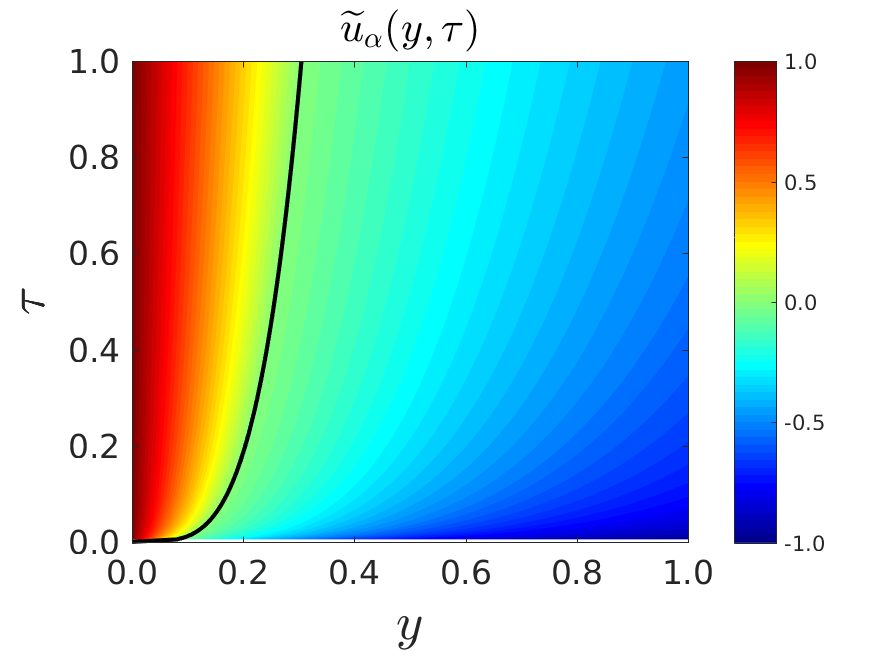} & 
\includegraphics[width=0.5\columnwidth]{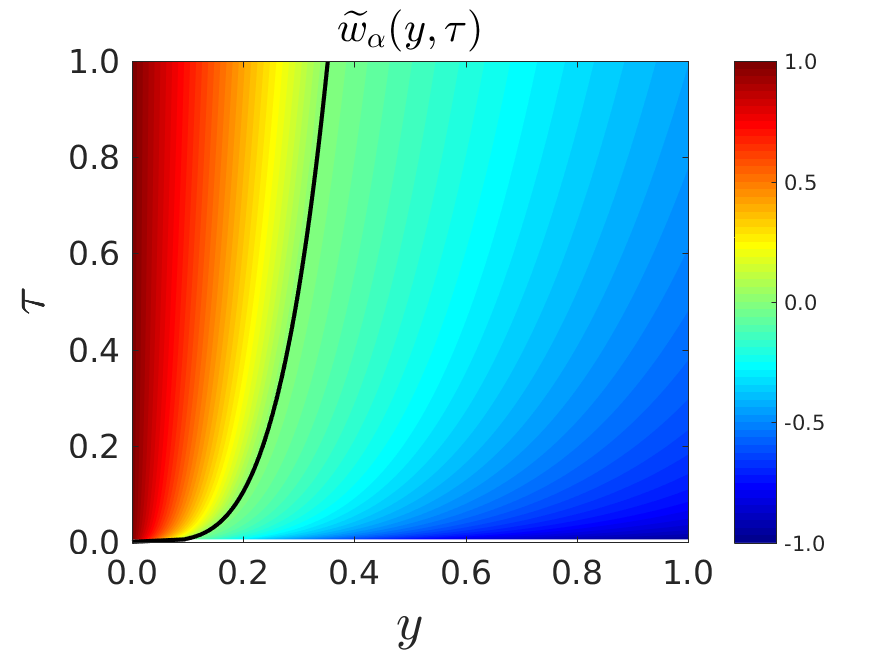} \\
$0.8$ & \includegraphics[width=0.5\columnwidth]{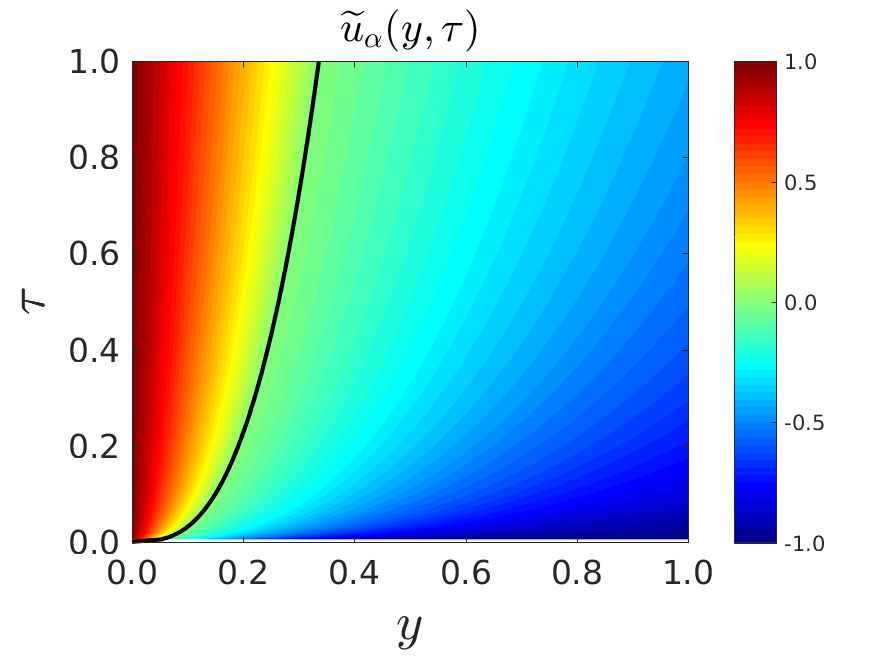} & 
\includegraphics[width=0.5\columnwidth]{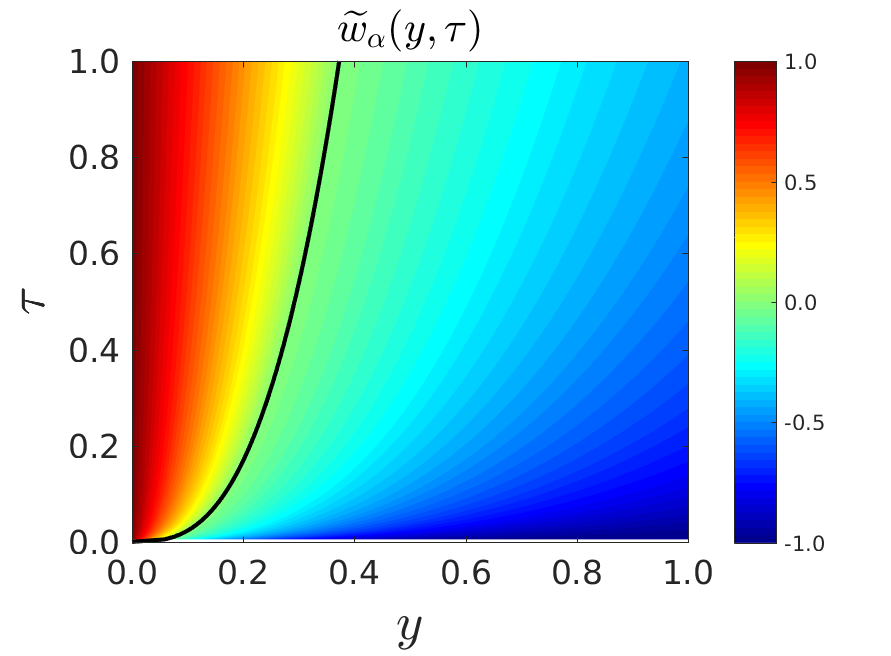} \\
$0.9$ & \includegraphics[width=0.5\columnwidth]{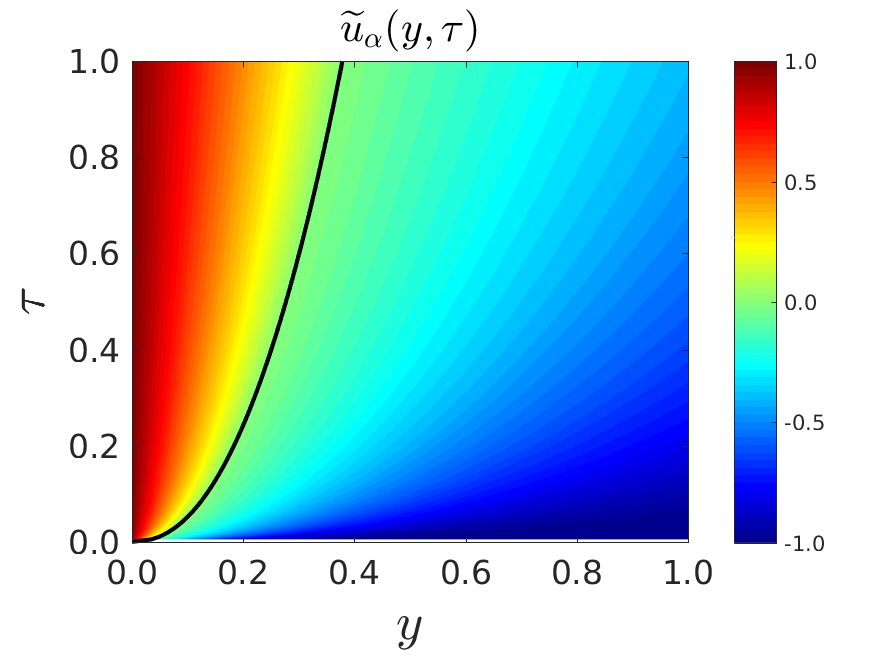} & 
\includegraphics[width=0.5\columnwidth]{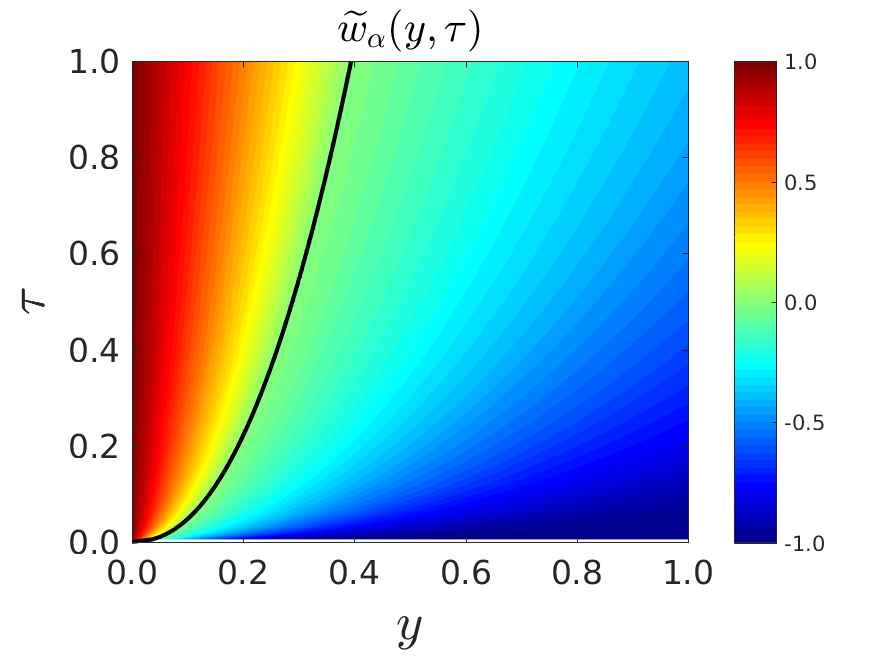}\\
\hline
\end{tabular}
\end{figure}

\begin {figure}[h]
\label{C_vs_RL_test3}
\caption {Caputos's approach Solutions Vs. Riemann-Liouville's aproach Solutions for Test 3} 
\centering
\begin{tabular}{|m{0.5cm}|m{7cm}| m{7cm}|} 

\hline
 $\alpha$ & \textbf{Caputo Stefan-Like Pb. } & \textbf{Riemann-Liouville  Stefan-Like Pb. } \\
\hline   
\hline
$0.7$ & \includegraphics[width=0.5\columnwidth]{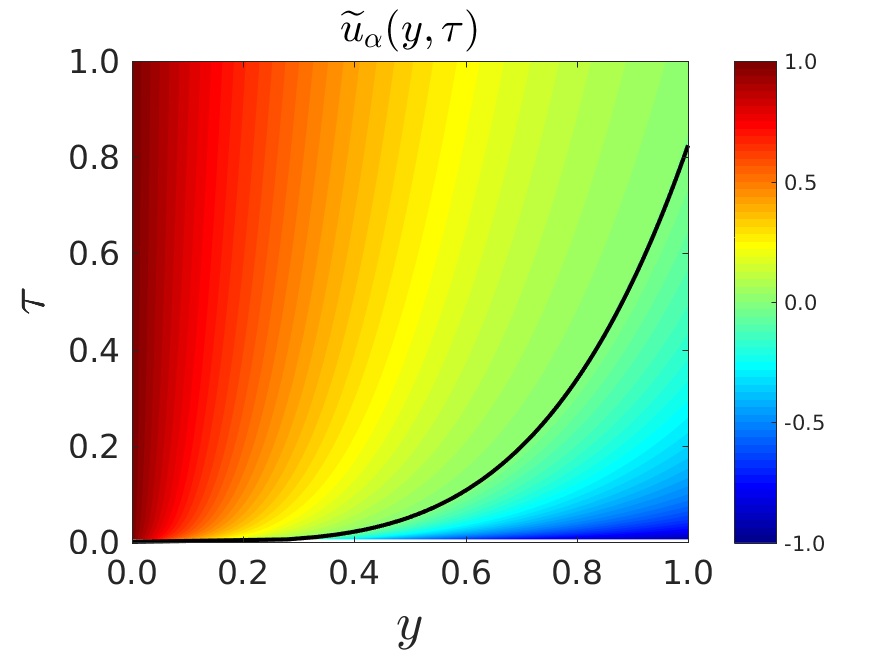} & 
\includegraphics[width=0.5\columnwidth]{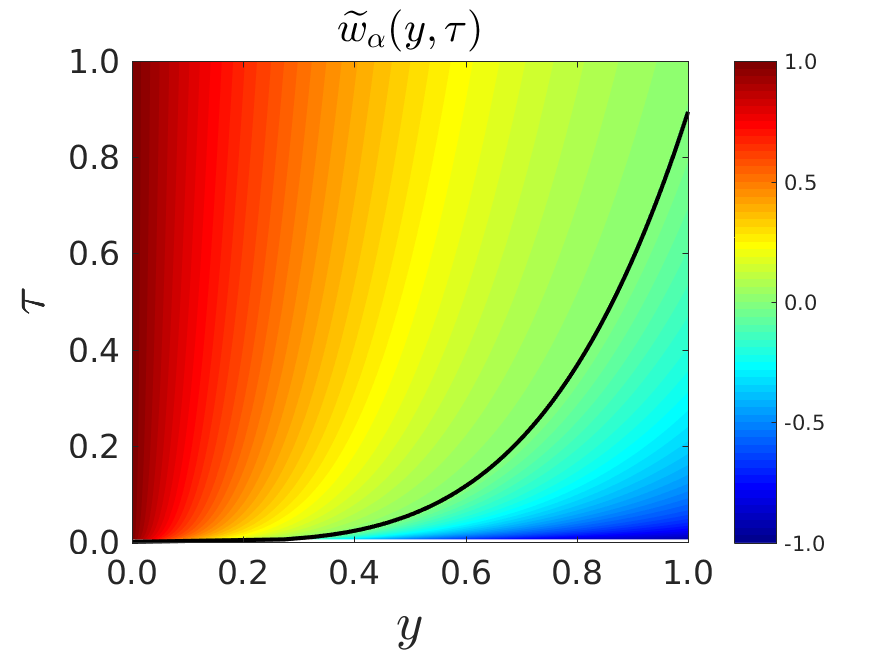} \\
$0.8$ & \includegraphics[width=0.5\columnwidth]{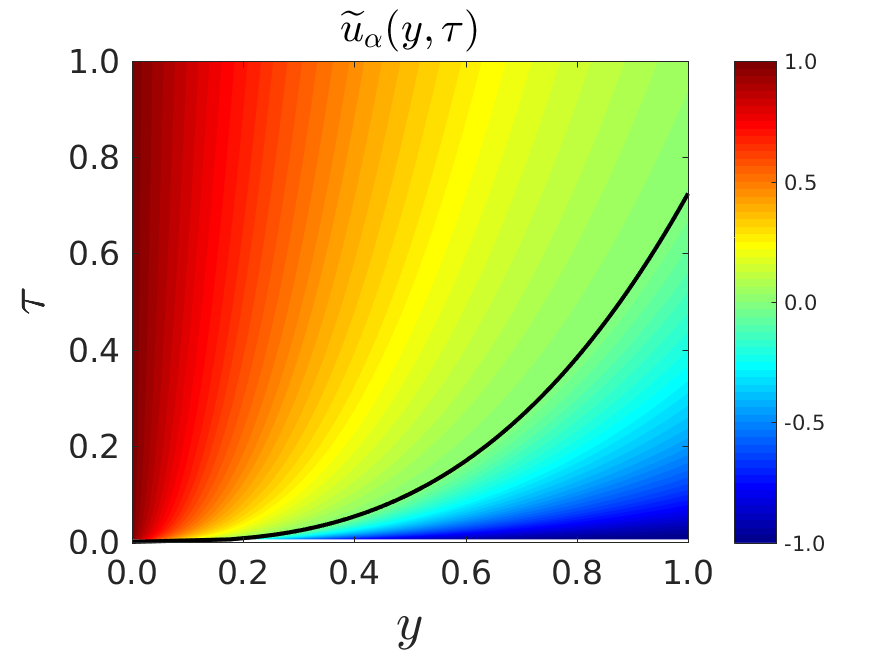} & 
\includegraphics[width=0.5\columnwidth]{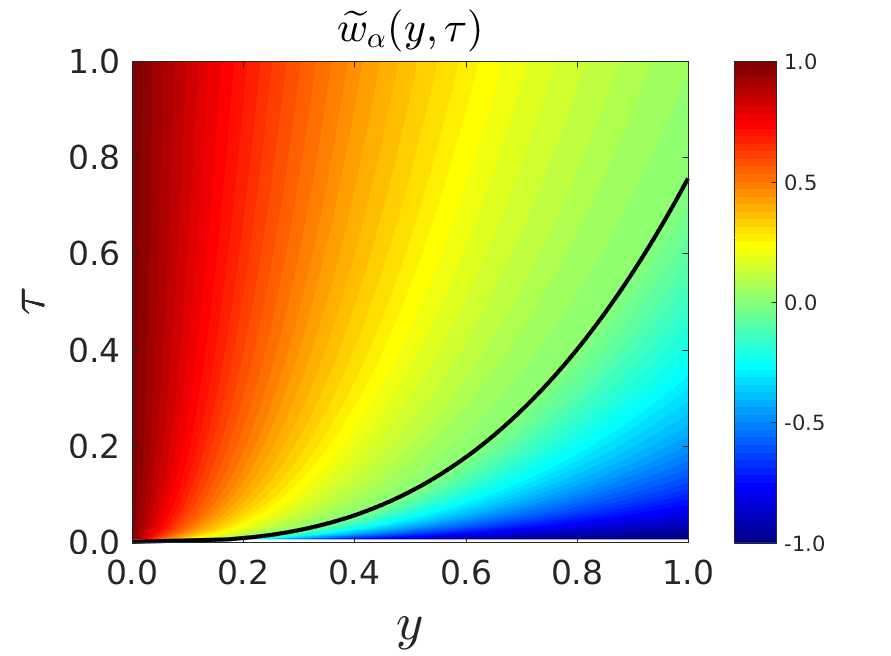} \\
$0.9$ & \includegraphics[width=0.5\columnwidth]{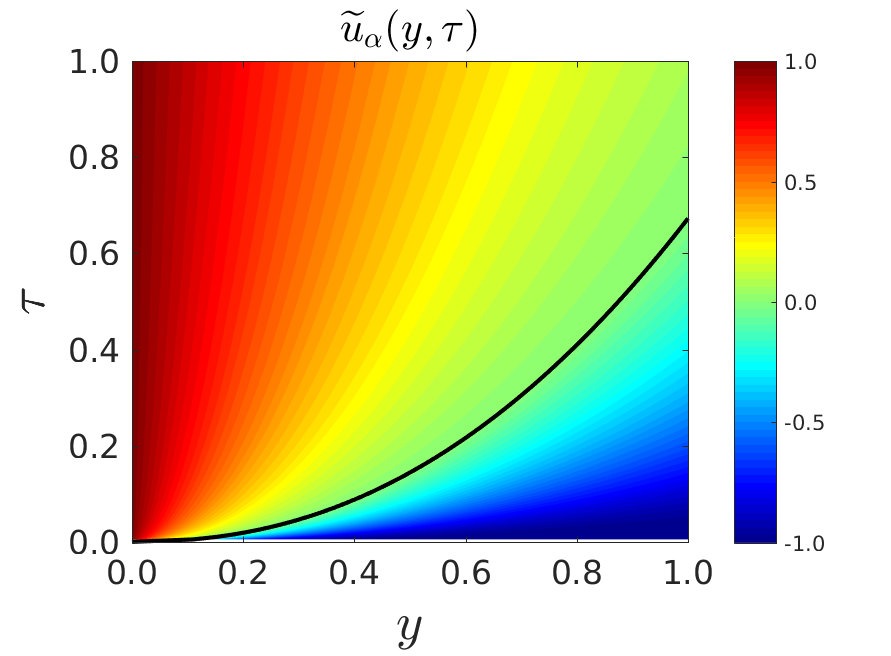} & 
\includegraphics[width=0.5\columnwidth]{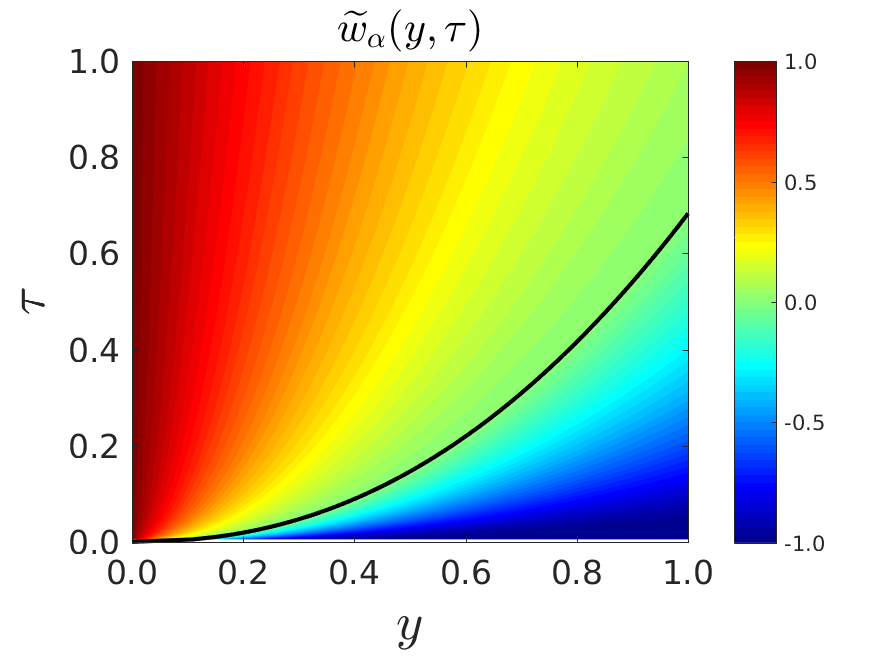}\\
\hline
\end{tabular}
\end{figure}

\newpage 

\section{Conclusion}
  We have presented two different fractional  two-phase Stefan-like problems  for the Riemann-Liouville and Caputo derivatives  of order $\al\in (0,1)$ with the particularity that, if the parameter $\al=1$ is replaced in both problems, we recover the same classical Stefan problem. In both cases, explicit solutions in terms of self-similar variables where given. 
It was interesting to see that, the role of the different  ``fractional Stefan conditions'' associated to each problem was decisive to conclude that the solutions obtained where different. Also, as it was expected, we have proved that the two different solutions converge to the same triple of limits functions when $\al$ tends to 1, where this ``limit solution'' is the classical solution to the classical Stefan problem mentioned before.

\section{Acknowledgements}

\noindent The present work has been sponsored by the Projects PIP N$^\circ$ 0275 from CONICET--Univ. Austral, ANPCyT PICTO Austral 2016 N$^\circ 0090$, Austral N$^\circ 006-$INV$00020$ (Rosario, Argentina) and European Unions Horizon 2020 research and innovation programme under the Marie Sklodowska-Curie Grant Agreement N$^\circ$ 823731 CONMECH.


\begin{thebibliography}{10}

\bibitem{Alexiades}
V.~Alexiades and A.~D. Solomon.
\newblock {\em Mathematical {M}odelling of {M}elting and {F}reezing
  {P}rocesses}.
\newblock Hemisphere, 1993.

\bibitem{BaFr:2005}
D.~S. Banks and C.~Fradin.
\newblock Anomalous diffusion of proteins due to molecular crowding.
\newblock {\em Biophysical Journal}, 89(5):2960 -- 2971, 2005.

\bibitem{GeVoMiDa:2013}
M.~Miksis C.~Gruber, C.~Vogl and S.~Davis.
\newblock Anomalous diffusion models in the presence of a moving interface.
\newblock {\em Interfaces and {F}ree {B}oundaries}, 15:181--202, 2013.

\bibitem{Cannon}
J.~R. Cannon.
\newblock {\em The {O}ne--{D}imensional {H}eat {E}quation}.
\newblock Cambridge University Press, 1984.

\bibitem{Crank}
J.~Crank.
\newblock {\em Free and {M}oving {B}oundary {P}roblems}.
\newblock Clarendon Press, 1984.

\bibitem{Diethelm}
K.~Diethelm.
\newblock {\em The {A}nalysis of {F}ractional {D}ifferential {E}quations: An
  application oriented exposition using differential operators of {C}aputo
  type}.
\newblock Springer Science \& Business Media, 2010.

\bibitem{Fujita:1989}
Y.~Fujita.
\newblock Integrodifferential equations which interpolates the heat equation
  and a wave equation.
\newblock {\em Osaka Journal of Mathematics}, 27:309--321, 1990.

\bibitem{GKS}
D.~N. Gerasimov, V.~A. Kondratieva, and O.~A. Sinkevich.
\newblock An anomalous non--self--similar infiltration and fractional diffusion
  equation.
\newblock {\em Physica D: Nonlinear Phenomena}, 239(16):1593--1597, 2010.

\bibitem{GoLuMa:2000}
Rudolf Gorenflo, Yuri Luchko, and Francesco Mainardi.
\newblock Wright functions as scale-invariant solutions of the diffusion-wave
  equation.
\newblock {\em Journal of Computational and Applied Mathematics},
  118(1):175--191, 2000.

\bibitem{JiMi:2009}
L.~Junyi and X.~Mingyu.
\newblock Some exact solutions to stefan problems with fractional differential
  equations.
\newblock {\em Journal of Mathematical Analysis and Applications},
  351:536--542, 2009.

\bibitem{FM-Libro}
F.~Mainardi.
\newblock {\em Fractional {C}alculus and {W}aves in {L}inear
  {V}iscoelasticity}.
\newblock Imperial Collage Press, 2010.

\bibitem{MK:2000}
R.~Metzler and J.~Klafter.
\newblock The random walk's guide to anomalous diffusion: a fractional dynamics
  approach.
\newblock {\em Physics reports}, 339:1--77, 2000.

\bibitem{Pa:2013}
G.~Pagnini.
\newblock The {M}-{W}right function as a generalization of the gaussian density
  for fractional diffusion processes.
\newblock {\em Fractional Calculus and Applied Analysis}, 16(2):436--453, 2013.

\bibitem{Podlubny}
I.~Podlubny.
\newblock {\em Fractional {D}ifferential {E}quations}.
\newblock Vol. 198 of Mathematics in Science and Engineering, Academic Press,
  1999.

\bibitem{Povstenko}
Y.~Povstenko.
\newblock {\em Linear {F}ractional {D}iffusion--wave {E}quation for
  {S}cientists and {E}ngineers}.
\newblock Springer, 2015.

\bibitem{Pskhu-Libro}
A.~V. Pskhu.
\newblock {\em Partial {D}ifferential {E}quations of {F}ractional {O}rder (in
  Russian)}.
\newblock Nauka, Moscow, 2005.

\bibitem{RoBoTa:2018}
S.~Roscani, J.~Bollati, and D.~Tarzia.
\newblock A new mathematical formulation for a {P}hase {C}hange {P}roblem with
  a memory flux.
\newblock {\em Chaos, Solitons and Fractals}, 116:340--347, 2018.

\bibitem{RoSa:2013}
S.~Roscani and E.~Santillan~Marcus.
\newblock Two equivalen {S}tefan's problems for the time--fractional diffusion
  equation.
\newblock {\em Fractional Calculus $\&$ Applied Analysis}, 16(4):802--815,
  2013.

\bibitem{RoTa:2014}
S.~Roscani and D.~Tarzia.
\newblock A generalized {N}eumann solution for the two--phase fractional
  {L}am\'e--{C}lapeyron--{S}tefan problem.
\newblock {\em Advances in Mathematical Sciences and Applications},
  24(2):237--249, 2014.

\bibitem{RoTa:2017-TwoDifferent}
S.~Roscani and D.~Tarzia.
\newblock Two different fractional {S}tefan problems which are convergent to
  the same classical {S}tefan problem.
\newblock {\em Mathematical Methods in the Applied Sciences}, 41(6):6842--6850,
  2018.

\bibitem{Samko}
S.~G. Samko, A.~A. Kilbas, and O.~I. Marichev.
\newblock {\em Fractional {I}ntegrals and {D}erivatives--{T}heory and
  {A}pplications}.
\newblock Gordon and Breach, 1993.

\bibitem{Sax:1994}
M.J. Saxton.
\newblock Anomalous diffusion due to obstacles: a monte carlo study.
\newblock {\em Biophysical Journal}, 66(2, Part 1):394 -- 401, 1994.

\bibitem{Tar:1981}
D.~A. Tarzia.
\newblock An inequality for the coeficient $\sigma$ of the free boundary
  $s(t)=2\sigma \sqrt{t}$ of the {N}eumann solution for the two-phase {S}tefan
  problem.
\newblock {\em Quart. Appl. Math.}, 39:491--497, 1981.

\bibitem{Tarzia:biblio}
D.~A. Tarzia.
\newblock A bibliography on moving--free boundary problems for the heat
  diffusion equation. the {S}tefan and related problems.
\newblock {\em MAT--Serie A}, 2:1--297, 2000.

\bibitem{Tarzia}
D.~A. Tarzia.
\newblock {\em Explicit and Approximated Solutions for Heat and Mass Transfer
  Problems with a Moving Interface}, chapter 20, {A}dvanced {T}opics in {M}ass
  {T}ransfer, pages 439--484.
\newblock Prof. Mohamed El-Amin (Ed.), Intech, Rijeka, 2011.

\bibitem{Voller:2014}
V.~R. Voller.
\newblock Fractional {S}tefan problems.
\newblock {\em International Journal of Heat and Mass Transfer}, 74:269--277,
  2014.

\bibitem{VoFaGa:2013}
V.~R. Voller, F.~Falcini, and R.~Garra.
\newblock Fractional {S}tefan problems exhibing lumped and distributed
  latent--heat memory effects.
\newblock {\em Physical Review E}, 87:042401, 2013.

\bibitem{Weber:1901}
H.~Weber.
\newblock {\em Die {P}artiellen {D}ifferential-{G}leichungen der
  {M}athematischen {P}hysik}.
\newblock Druck und Verlag von Friedrich Vieweg und Sohn, 1901.

\bibitem{Wr1:1934}
E.~M. Wright.
\newblock The assymptotic expansion of the generalized bessel funtion.
\newblock {\em Proceedings of the London mathematical society (2)},
  38:257--270, 1934.

\bibitem{Wr2:1940}
E.~M. Wright.
\newblock The generalized {B}essel function of order greater than one.
\newblock {\em The Quarterly Journal of Mathematics}, Ser. 11:36--48, 1940.

\end{thebibliography}
\end{document}